\documentclass[12pt,a4j]{paper}

\usepackage{amsmath,amssymb}
\usepackage{amsthm}
\usepackage{amsrefs}
\usepackage{cancel}
\usepackage{color}
\usepackage[legacycolonsymbols]{mathtools}
\usepackage{tikz}
\usepackage{tikz-3dplot}
\usepackage[compat=1.1.0]{tikz-feynhand}
\usepackage{empheq}

\numberwithin{equation}{section}
\newtheorem{thm}{Theorem}[section]

\newtheorem{dfn}[thm]{Definition}
\theoremstyle{definition}
\newtheorem{ex}[thm]{Example}
\newtheorem{rem}[thm]{Remark}
\newcommand{\del}{\partial}

\newcommand{\df}{\coloneqq}
\newcommand{\beq}{\begin{equation}}
\newcommand{\eeq}{\end{equation}}
\newcommand{\bsp}{\begin{split}}
\newcommand{\esp}{\end{split}}
\newcommand{\redcancel}[1]{\textcolor{red}{\cancel{\textcolor{black}{#1}}}}

\newcommand{\redunderline}[1]{\textcolor{red}{\underline{\textcolor{black}{#1}}}}
\newcommand{\blueunderline}[1]{\textcolor{blue}{\underline{\textcolor{black}{#1}}}}

\begin{document}

\begin{center}
{\Large  Topological invariants of 3-dimensional manifold with boundary by using crossed module}
\\
Tommy Shu
\end{center}

\section{introduction}
Crossed module is invited by J.H.C. Whitehead(\cite{Whi49}) in early 20 century. This is induced for algebra model of 2-type homotopy(i.e connected spaces where there are no more than 2 degree homotopy group). Crossed module, which is define by two group and some relation between them, is known that it have a relation with 2-group.

By using crossed module, we can construct invariant of closed  3-dimensional manifold and 4-dimensional manifold which proof is witten by (\cite{FHE08}). Roughly, this invariant is counting correct coloer over the trianguler of closed manifold. 

I find that this invariant can be used to compact 3-dimensional manifold with boundary. So in this paper, I will proof that this invariant can  be used to compact 3-dimensional manifold with boundary and show some result of using this invariant to  manifold with boundary(inclue  knot complement).
 
 \section{Crossed module}

Crossed module use two group, homeomorphism, and left action. 

\begin{dfn}
Let $H,G$ be a group, $\del:H\to G$ be a homeomorhism, $\rhd$ as a left action of $G$ on $G$. and  $\rhd$ as a left action of $G$ on $H$(we write both left action as a same sympol because most of the time it is obvious). We call $(H\stackrel{\del}{\to}G,\rhd)$ crossed module if it satisfies the following three properties:
\begin{enumerate}
	 \item
	 all $g_1,g_2 \in G$ satisfy
	 \[
	 g_1 \rhd g_2 = g_1 g_2 g_1^{-1};
	 \]
	 \item
	 all $g \in G$ and  for all $h \in H$ satisfy
	 \[
	 \del (g \rhd h) = g\del(h)g^{-1};
	 \]
	 and
	 \item if a homeomorhism $f_g:H \to H$ is defined by $f_g(h)\df g\rhd h $ for $h\in H$, then
	 $f_g$ is a group isomorphism, for all $g \in G$.
\end{enumerate}
\label{df,cm}
\end{dfn}

 We will give a most simple example of crossed module.
 \begin{ex}
 Let $H$ be a finite group, and let $G$ be a {\rm Aut}$(H)$$\df\{f:H\to H|$$f$ is group isomorphism\}. Then we have a crossed module $(H\stackrel{\del}{\to}G,\rhd)$  if $\del$ and $ \rhd$ are defined as follows:
\begin{enumerate}
	 \item
	  for all $g_1,g_2 \in G$ we define
	 \[
	 g_1 \rhd g_2 \df g_1 g_2 g_1^{-1};
	 \]
	 \item
	 for all $g \in G$ and all $h \in H$  define
	 \[
	 g \rhd h \df g(h);
	 \]
	 and
	 \item for all $h,h'\in H$ define
	 \[
	 \del(h')(h)\df h'hh'^{-1}.
	 \]
\end{enumerate}
 \end{ex}
 
 In the next section we use this crossed module to construct an invariant.
 
  \section{Invariant of closed 3-dimensional manifold}
 To know more about this section you can see   (\cite{TM22}) and (\cite{FHE08}). The way to construct an invariant with crossed modole gives by next theorem.
 
\begin{thm}
Let M be a closed 3-dimensional manifold and let $(H\stackrel{\del}{\to}G,\rhd )$ be a crossed module but both H and G is finite group. K is a triangular of M ,$K_i$ is a set of all i-simplex in K , and there are total order in $K_0$ .Than next formuler give an invariant of M and it don't depend on choice of K and total order in $K_0$:
\begin{equation}
\begin{split}
Z(M)=&|G|^{-|K_0|+|K_1|-|K_2|}|H|^{|K_0|-|K_1|+|K_2|-|K_3|}\\
&\times\left(\prod_{(jk)\in K_1}{1\over |G|}\sum_{g_{jk}\in G}\right)\left(\prod_{(jkl)\in K_2}{1\over |H|}\sum_{h_{jkl}\in H}\right)\\
&\times\left(\prod_{(jkl)\in K_2} \delta_G(g_{jkl})\right) \left(\prod_{(jklm)\in K_3} \delta_H(h_{jklm})\right)
\end{split}
\eeq
where 
\beq
g_{jkl}\coloneqq \del(h_{jkl})g_{kl}g_{jk}g_{jl}^{-1}
\label{gjkl}
\eeq
\beq
h_{jklm}\coloneqq h_{jlm}(g_{lm}\rhd{ h_{jkl}})h_{klm}^{-1}h_{jkm}^{-1},
\label{hjklm}
\eeq
and, $\delta_X:X\to \mathbb{Z}$ for finite group X which define by
\begin{equation}
  \delta_X(x)=
  \begin{cases}
    |X|  & \text {if $x=e_X$,} \\
    0              &    \text{others.} .
  \end{cases}
\end{equation}
However, the label of $g_{jk}, h_{jkl}$ will always be $j<k$ and $j<k<l$ of total order in $K_0$.
\label{closed 3mfd thm}
\end{thm}
  
\begin{figure}[htbp]
\begin{center}
\begin{tikzpicture}
\begin{feynhand}
\propag[fer] (2,2) to (0,0);
\propag[fer] (3,0) to (2,2);
\propag[fer] (3,0) to (0,0);
\draw(0,0)node[left]{1};
\draw(2,2)node[above]{2};
\draw(3,0)node[right]{3};
\draw(1,1)node[left]{$g_{12}$};
\draw(2.5,1)node[right]{$g_{23}$};
\draw(1.5,0)node[below]{$g_{13}=\del(h_{123})g_{23}g_{12}$};
\draw(1.7,0.6)node{$h_{123}$};
\end{feynhand}
\end{tikzpicture}
\end{center}
\caption{Color of $g_{123}$}
\label{2color}
\end{figure}
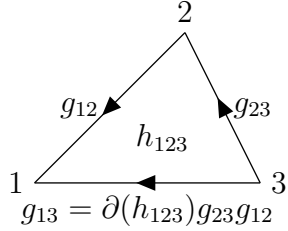

Proof can be finded in \cite{FHE08}. What this invariant doing is counting the correct color. We color $g_{jk}\in G$ for evey edge $(jk)\in K_1$  and $h_{jkl}\in H$ for every face $(jkl)\in K_1$. The geometric meaning of $g_{jkl}$ become figure \ref{2color}.

Since there are total order in $K_0$, there are  no cyclic order for 2-simplex. For the order of the each edge which induced by total order in $K_0$, as you can see in  figure \ref{2color}, $g_{13}$ can be written by other two colors of edges and one color of face. 

\begin{figure}
\begin{center}
\includegraphics[width=50mm]{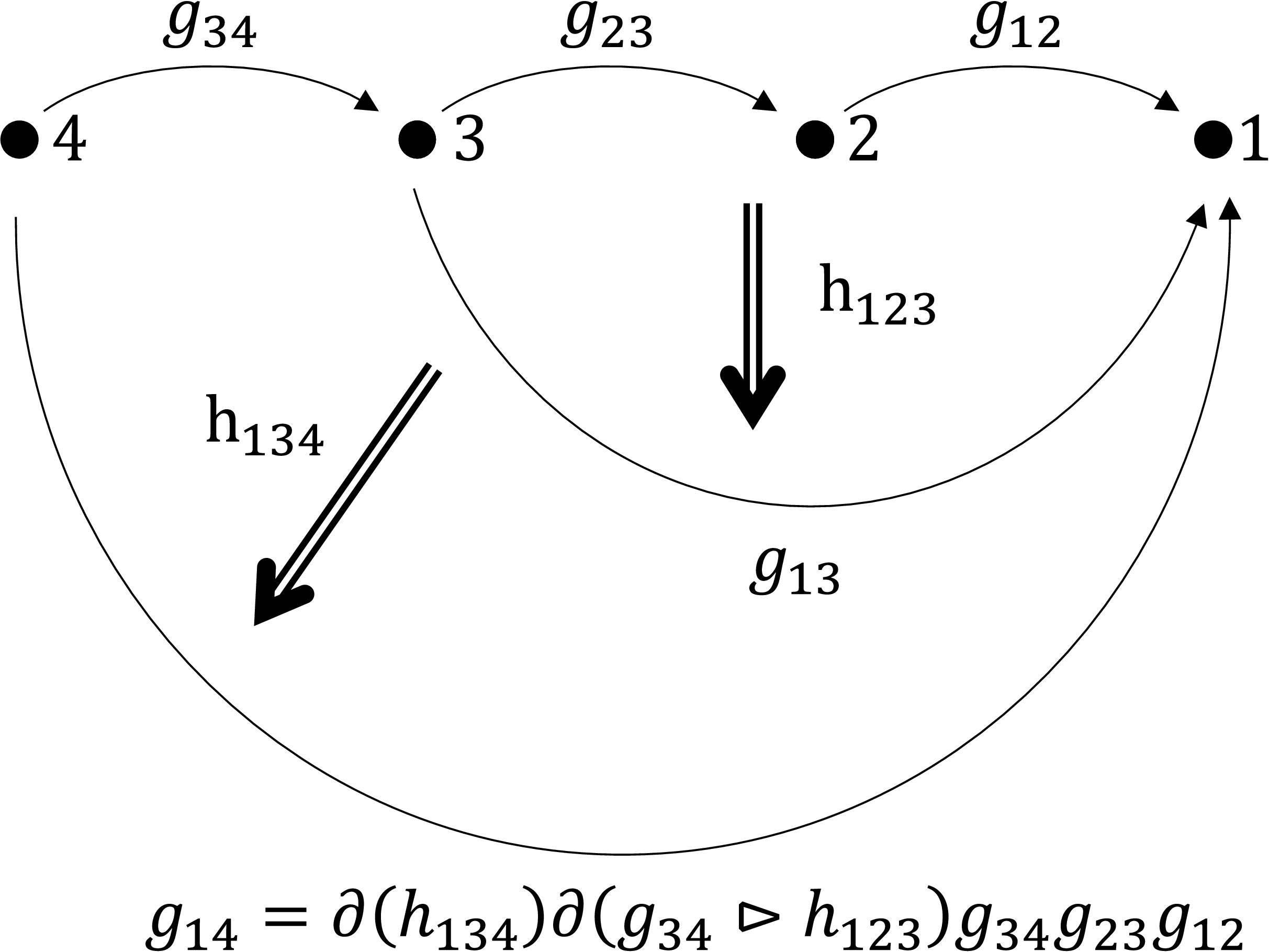}
\caption{Using $g_{12}, g_{23}, g_{34}, h_{123} ,h_{134}$ }
\label{1w}
\end{center}
\end{figure}

\begin{figure}
\begin{center}
\includegraphics[width=50mm]{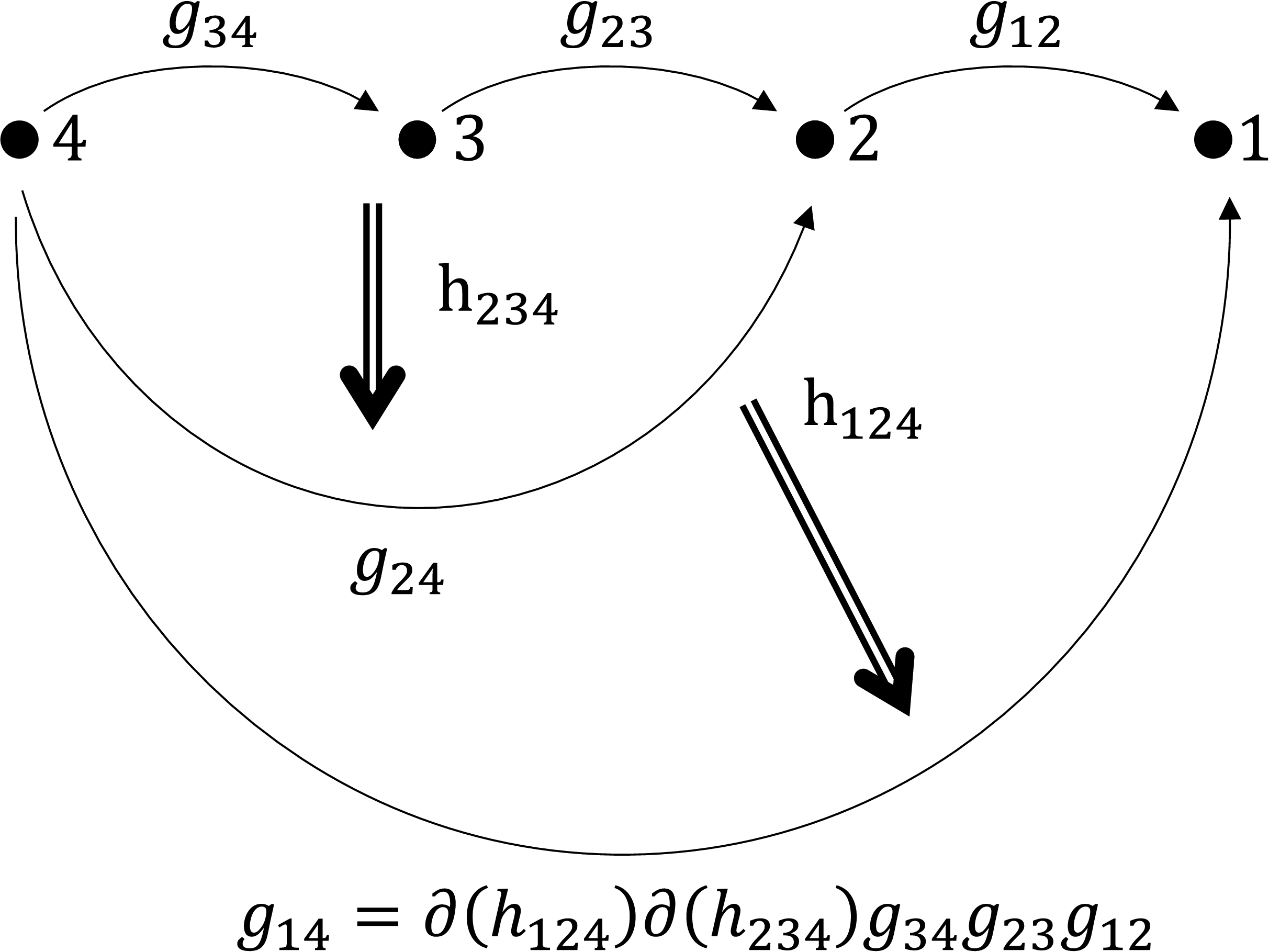}
\caption{Using $g_{12}, g_{23}, g_{34}, h_{234} ,h_{124}$}
\label{2w}
\end{center}
\end{figure}

We can write $g_{14}$ in two ways in which one is used $g_{12}, g_{23}, g_{34}, h_{123} ,h_{134}$(diagram \ref{1w}) and the other is used $g_{12}, g_{23}, g_{34}, h_{234} ,h_{124}$
(see Figure~\ref{2w}).

\begin{figure}
\begin{center}
\includegraphics[width=50mm]{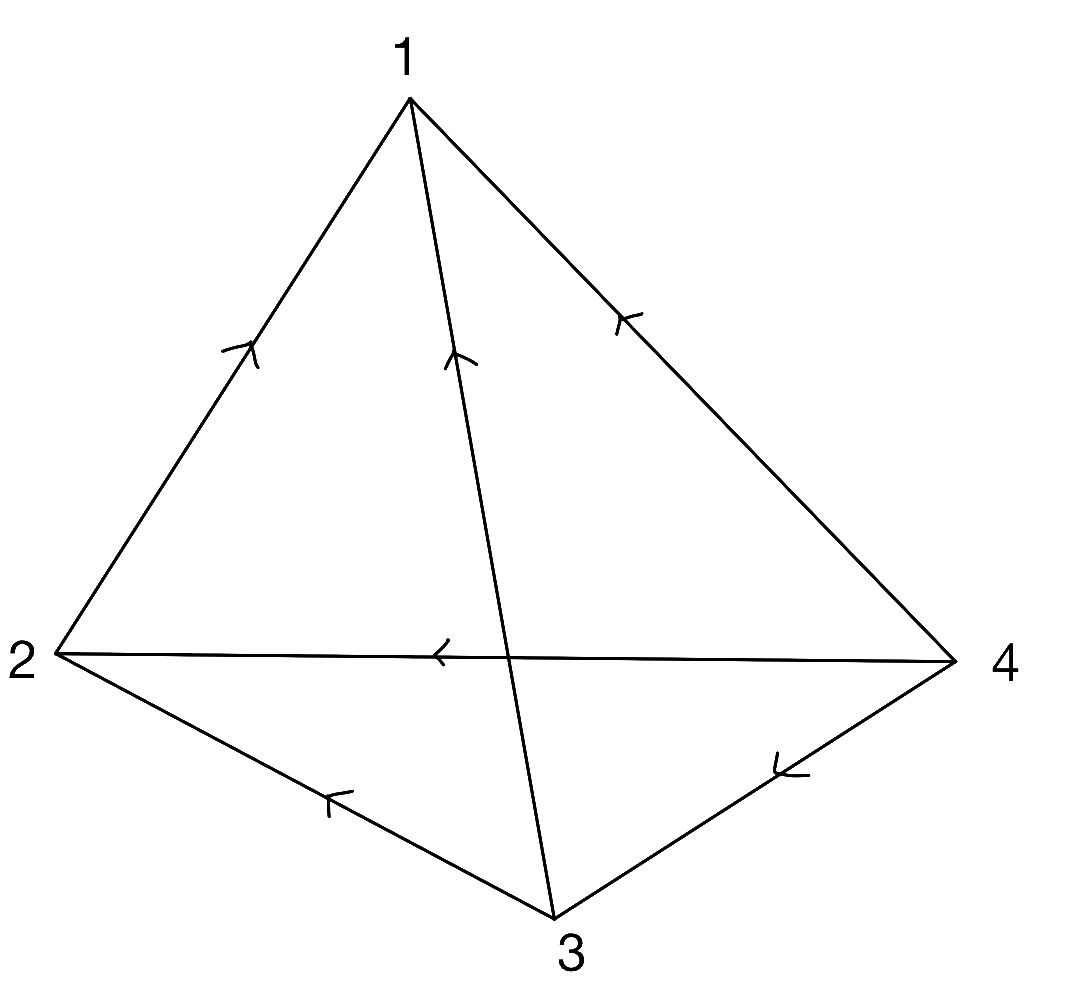}
\caption{}
\label{mhjklm}
\end{center}
\end{figure}

This two formula must be equation so we get next equation.

\beq
\del\left(h_{134}(g_{34}\rhd h_{123})\right)=\del(h_{124}h_{234}).
\label{3tanntainocomp}
\eeq

When formula \ref{hjklm} is equal to identity element, equation \ref{3tanntainocomp} will hold true. 

\section{Main Theorem}

We extend the theorem \ref{closed 3mfd thm} to not only closed 3-dimensional manifold but also compact 3-dimensional manifold with boundary.

\begin{thm}
Let M be a \redunderline{compact 3-dimensional manifold with boundary} and let $(H\stackrel{\del}{\to}G,\rhd )$ be a crossed module but both H and G is finite group. K is a triangular of M ,$K_i$ is a set of all i-simplex in K , and there are total order in $K_0$ .Than next formuler give an invariant of M and it don't depend on choice of K and total order in $K_0$:
\begin{equation}
\begin{split}
Z(M)=&|G|^{-|K_0|+|K_1|-|K_2|}|H|^{|K_0|-|K_1|+|K_2|-|K_3|}\\
&\times\left(\prod_{(jk)\in K_1}{1\over |G|}\sum_{g_{jk}\in G}\right)\left(\prod_{(jkl)\in K_2}{1\over |H|}\sum_{h_{jkl}\in H}\right)\\
&\times\left(\prod_{(jkl)\in K_2} \delta_G(g_{jkl})\right) \left(\prod_{(jklm)\in K_3} \delta_H(h_{jklm})\right)
\end{split}
\label{ZM}
\eeq
where 
\beq
g_{jkl}\coloneqq \del(h_{jkl})g_{kl}g_{jk}g_{jl}^{-1}
\label{gjkl}
\eeq
\beq
h_{jklm}\coloneqq h_{jlm}(g_{lm}\rhd{ h_{jkl}})h_{klm}^{-1}h_{jkm}^{-1},
\label{hjklm}
\eeq
and, $\delta_X:X\to \mathbb{Z}$ for finite group X which define by
\begin{equation}
  \delta_X(x)=
  \begin{cases}
    |X|  & \text {if $x=e_X$,} \\
    0              &    \text{others.} .
  \end{cases}
\end{equation}
However, the label of $g_{jk}, h_{jkl}$ will always be $j<k$ and $j<k<l$ of total order in $K_0$.
\label{3mfd_w_b thm}
\end{thm}

What this theorem want to say is it doesn't depend the triangular of its boundary. Before moving to the proof, I will give two theoremes that is used to proof the theorem \ref{3mfd_w_b thm}.

\begin{figure}[htbp]
\begin{center}
\tdplotsetmaincoords{60}{0}
\begin{tikzpicture}[tdplot_main_coords,scale=1]
\coordinate (A) at (2,0,0);
\coordinate (B) at ({2*cos(120)},{2*sin(120)},0);
\coordinate (C) at ({2*cos(-120)},{2*sin(-120)},0);
\coordinate (D) at (0,0,3);
\coordinate (AA) at (0,0,-3);
\draw(A)--(B)--(C)--cycle;
\draw(A)--(C)--(D)--cycle;
\draw(A)--(B)--(D)--cycle;
\draw(D)--(B)--(C)--cycle;
\draw(AA)--(C);
\draw(AA)--(B);
\draw(AA)--(A);
\coordinate (E) at (2+6,0,0);
\coordinate (F) at ({2*cos(120)+6},{2*sin(120)},0);
\coordinate (G) at ({2*cos(-120)+6},{2*sin(-120)},0);
\coordinate (H) at (6,0,3);
\coordinate (EE) at (6,0,-3);
\draw(E)--(F)--(G)--cycle;
\draw(E)--(G)--(H)--cycle;
\draw(E)--(F)--(H)--cycle;
\draw(F)--(G)--(H)--cycle;
\draw(EE)--(E);
\draw(EE)--(F);
\draw(EE)--(G);
\draw(EE)--(H);
\draw[<->] (3,0,0)--(4,0.0);
\draw (7/2,0,0)node[above]{(2-3)};
\end{tikzpicture}
\end{center}
\caption{(2-3)move}
\label{fig:(2-3)}
\end{figure}
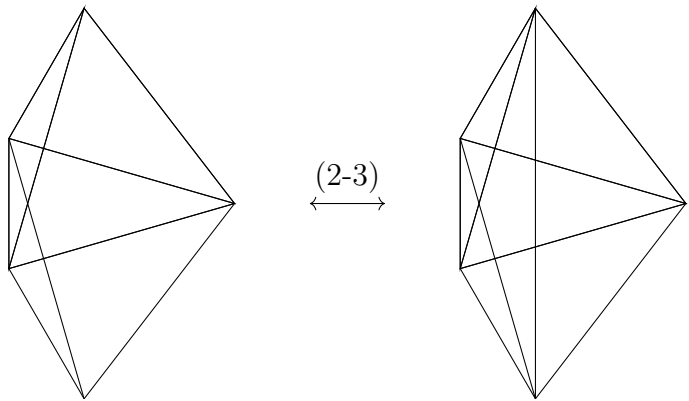

\begin{thm}(Turaev-Viro)
\label{th:tv}

Let $P ,Q$ be the n-dimensional  Combinatorial manifold and $K$and $L$ is a triangular of $P$ and $Q$ which satisfy $ P = |K|, Q = |L|$. In this case next two properties is equivalent:
\begin{enumerate}
\item There exists a piecewise linear homeomorphism $f:P\to Q$ where $f|_{\del P}:\del P \to \del Q$ is a simplex homeomorphism from $\del K$ to $\del L$.
\item $L$ can be obteined by using finite Alexander moves from $K$.
\\

\end{enumerate}
\end{thm}
In the case where 
$n=3$, it has been showen that internal Alexander moves can be realized through a finite sequence of 3-dimensional Pachner moves.(\cite{iop96}:Theorem 4.6).
 Figure \ref{fig:(1-4)}, \ref{fig:(2-3)} is 3-dimensional Pachner moves.

\begin{figure}[htbp]
\begin{center}
\tdplotsetmaincoords{60}{0}
\begin{tikzpicture}[tdplot_main_coords,scale=1]
\coordinate (A) at (2,0,0);
\coordinate (B) at ({2*cos(120)},{2*sin(120)},0);
\coordinate (C) at ({2*cos(-120)},{2*sin(-120)},0);
\coordinate (D) at (0,0,3);
\draw(A)--(B)--(C)--cycle;
\draw(A)--(C)--(D)--cycle;
\draw(A)--(B)--(D)--cycle;
\draw(D)--(B)--(C)--cycle;
\coordinate (E) at (2+6,0,0);
\coordinate (F) at ({2*cos(120)+6},{2*sin(120)},0);
\coordinate (G) at ({2*cos(-120)+6},{2*sin(-120)},0);
\coordinate (H) at (6,0,3);
\coordinate (EE) at ($0.25*(E)+0.25*(F)+0.25*(G)+0.3*(H)$);
\draw(E)--(F)--(G)--cycle;
\draw(E)--(G)--(H)--cycle;
\draw(E)--(F)--(H)--cycle;
\draw(F)--(G)--(H)--cycle;
\draw(EE)--(E);
\draw(EE)--(F);
\draw(EE)--(G);
\draw(EE)--(H);
\draw[<->] (3,1,0)--(4,1.0);
\draw (7/2,1,0)node[above]{(1-4)};
\end{tikzpicture}
\end{center}
\caption{(1-4)move}
\label{fig:(1-4)}
\end{figure}
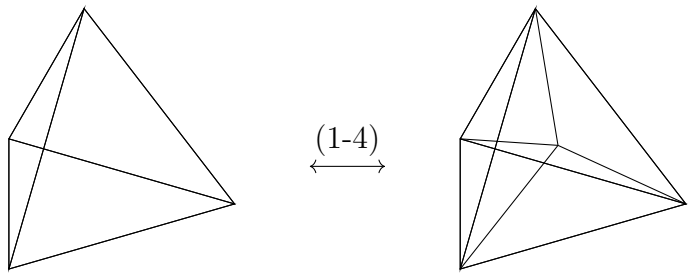

\begin{rem}
By \cite{moise}, every triangulated $3$-manifold $K$ (with boundary) is a combinatorial $3$-manifold (with boundary). $K$ is triangulated $n$-manifold, if K is a complex satisfying  that the space $M =|K|$ is an $n$-manifold. 
\end{rem}

\begin{thm}(\cite{mun}:Theorem 10.6)
\label{th:btri}

Let $M$ be a manifold with boundary. In this case, for any simplex subdivision of the boundary $\del M$ , there exists an extended triangular of $M$.The meaning of extended triangular is as follow: Let $J$  be a simplicial complex, and consider a homeomorphism $f:J \to \del M$ as a simplex subdivision of $\del M$. In this context, an extension of $f$ efers to the existence of a simplicial complex $L$ and a homeomorphism $g:L\to M$ with a simplex subdivision of $M$, such that $g^{-1}\circ f$ induces a linear isomorphism between a subcomplex of $L$ and$J$.
\end{thm}

\section{Proof of Main Theorem}

\begin{proof}

In the appendix A, we prove that formula \ref{ZM} does not depend on the choice of total order in $K_0$ and it is invariant in 3-dimension Pachner move. Since of this fact and theorem \ref{th:tv}, if there are two different triangular $K,L$ of $M$ which $\del K$ and $\del L$ of triangular of $\del M$ coincide, then $Z_K(M)=Z_L(M)$ which one is used K and the other used L. Only we have to prove is to show that we can transform two different triangular of $\del M$ one to another by preserving $Z(M)$. To transform the triangular of the boundary, we think about figure \ref{1-3mv}, \ref{2-2mv}.

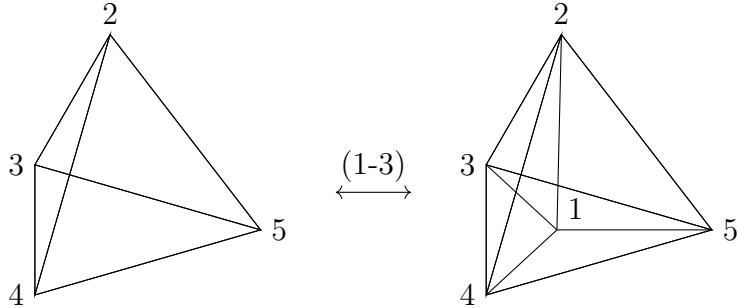
\begin{figure}[htbp]
\begin{center}
\tdplotsetmaincoords{60}{0}
\begin{tikzpicture}[tdplot_main_coords,scale=1]
\coordinate (A) at (2,0,0);
\coordinate (B) at ({2*cos(120)},{2*sin(120)},0);
\coordinate (C) at ({2*cos(-120)},{2*sin(-120)},0);
\coordinate (D) at (0,0,3);
\draw(A)--(B)--(C)--cycle;
\draw(A)--(C)--(D)--cycle;
\draw(A)--(B)--(D)--cycle;
\draw(D)--(B)--(C)--cycle;
\coordinate (E) at (2+6,0,0);
\coordinate (F) at ({2*cos(120)+6},{2*sin(120)},0);
\coordinate (G) at ({2*cos(-120)+6},{2*sin(-120)},0);
\coordinate (H) at (6,0,3);
\coordinate (EE) at ($0.33*(E)+0.33*(F)+0.33*(G)$);
\draw(E)--(F)--(G)--cycle;
\draw(E)--(G)--(H)--cycle;
\draw(E)--(F)--(H)--cycle;
\draw(F)--(G)--(H)--cycle;
\draw(EE)--(E);
\draw(EE)--(F);
\draw(EE)--(G);
\draw(EE)--(H);
\draw[<->] (3,1,0)--(4,1.0);
\draw(A)node[right]{5};
\draw(B)node[left]{3};
\draw(C)node[left]{4};
\draw(D)node[above]{2};
\draw(E)node[right]{5};
\draw(F)node[left]{3};
\draw(G)node[left]{4};
\draw(H)node[above]{2};
\draw(EE)node[above right]{1};
\draw (7/2,1,0)node[above]{(1-3)};
\end{tikzpicture}
\end{center}
\caption{3-dimension (1-3)move}
\label{1-3mv}
\end{figure}

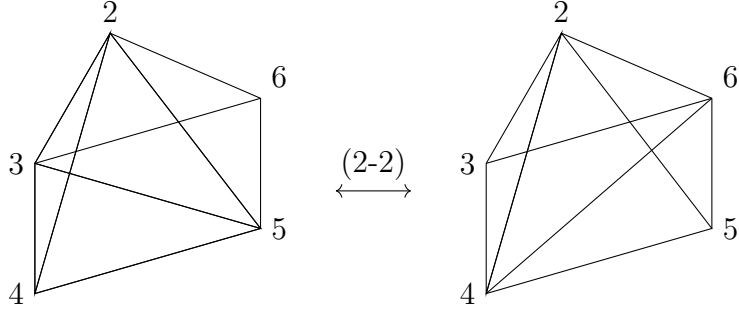
\begin{figure}[htbp]
\begin{center}
\tdplotsetmaincoords{60}{0}
\begin{tikzpicture}[tdplot_main_coords,scale=1]
\coordinate (A) at (2,0,0);
\coordinate (B) at ({2*cos(120)},{2*sin(120)},0);
\coordinate (C) at ({2*cos(-120)},{2*sin(-120)},0);
\coordinate (D) at (0,0,3);
\coordinate (AA) at ($(B)+(A)-(C)$);
\draw(A)--(B)--(C)--cycle;
\draw(A)--(C)--(D)--cycle;
\draw(A)--(B)--(D)--cycle;
\draw(D)--(B)--(C)--cycle;
\draw(AA)--(A);
\draw(AA)--(B);
\draw(AA)--(D);
\coordinate (E) at (2+6,0,0);
\coordinate (F) at ({2*cos(120)+6},{2*sin(120)},0);
\coordinate (G) at ({2*cos(-120)+6},{2*sin(-120)},0);
\coordinate (H) at (6,0,3);
\coordinate (EE) at ($(E)+(F)-(G)$);

\draw(E)--(G)--(H)--cycle;

\draw(F)--(G)--(H)--cycle;
\draw(EE)--(E);
\draw(EE)--(F);
\draw(EE)--(H);
\draw(EE)--(G);
\draw[<->] (3,1,0)--(4,1.0);
\draw(AA)node[above right]{6};
\draw(A)node[right]{5};
\draw(B)node[left]{3};
\draw(C)node[left]{4};
\draw(D)node[above]{2};
\draw(E)node[right]{5};
\draw(F)node[left]{3};
\draw(G)node[left]{4};
\draw(H)node[above]{2};
\draw(EE)node[above right]{6};
\draw (7/2,1,0)node[above]{(2-2)};
\end{tikzpicture}
\end{center}
\caption{3-dimension (2-2)move}
\label{2-2mv}
\end{figure}
However, there are conditions for performing this move. The (1-3) move in Figure \ref{1-3mv} is permissible when the 2-simplex |345| belongs to the boundary $\del M$ of manifold $M$,i.e. when $|345| \subset \del M$.
As for the (2-2) move in Figure \ref{2-2mv}, it can be executed when the union of 2-simplices |245| and |356| belongs to the boundary $\del M$,i.e. when $|245|\cup|356|\subset \del M$. The 3-dimension (1-3)move and (2-2)move preserve formula \ref{ZM} which you can see in appendix------.
 
Next we will prove that two different triangular of $M$ which we write $K,L$ can be transformed one to another by a sequuence of 3-dimension (1-3)move, (2-2) move, and pachner move. As privious stated, It is sufficient to show that  $\del K$ can transform to $\del L$ by a sequuence of 3-dimension (1-3)move, (2-2) move, and pachner move. $\del K$ and  $\del L$ is triangular of $\del M$ so we can be think as triangular of closed 2-dimensional manifold. If we forget about internal of $M$, $\del K$ can be transformed to $\del L$ by sequence of 2-dimension Pachner move which show in Figure \ref{2d13mv},\ref{2d22mv}.

\begin{figure}[htbp]
\begin{center}
\begin{tikzpicture}
\draw(0,0)--(2,2)--(3,0)--cycle;

\draw(6,0)--(8,2)--(9,0)--cycle;
\draw(8,2/3)--(6,0);
\draw(8,2/3)--(8,2);
\draw(8,2/3)--(9,0);
\draw[<->] (4,1)--(5,1);
\draw(4.5,1)node[above]{};
\end{tikzpicture}
\end{center}
\caption{2-dimension (1-3)move}
\label{2d13mv}
\end{figure}
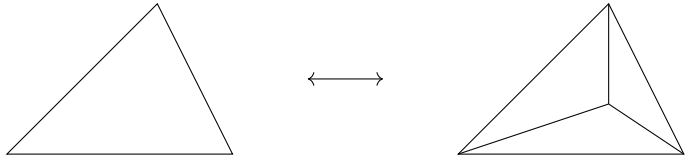
\begin{figure}[htbp]
\begin{center}
\begin{tikzpicture}
\draw(0,0)--(2,2)--(3,0)--cycle;
\draw(0,0)--(1,-2)--(3,0);
\draw(6,0)--(8,2)--(9,0);
\draw(6,0)--(7,-2)--(9,0);
\draw(8,2)--(7,-2);
\draw[<->] (4,0)--(5,0);
\draw(4.5,1)node[above]{};
\end{tikzpicture}
\end{center}
\caption{2-dimension (2-2)move}
\label{2d22mv}
\end{figure}
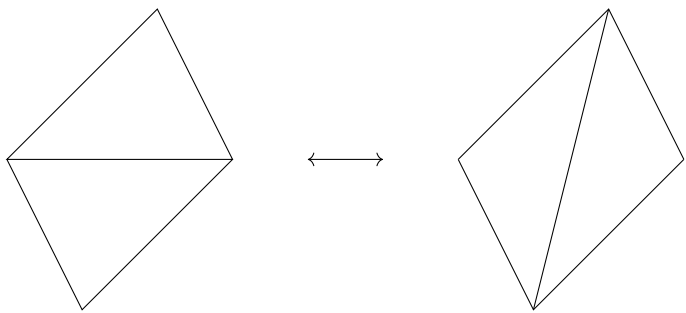

It is sufficient to prove that the sequence of two-dimensional Pachner moves, consisting of two-dimensional (1-3) moves and two-dimensional (2-2) moves, can be replaced with corresponding three-dimensional (1-3) moves and three-dimensional (2-2) moves.
\\

{(a) About (1-3)move}

It is trivial that we can replaced two-dimensional (1-3) moves from left to right with three-dimensional (1-3) moves from left to right. The reason why the transformation of a two-dimensional Pachner move (1-3) move from right to left can be replaced with a three-dimensional (1-3) move is as follows: For the portion where you want to transform the two-dimensional Pachner move (1-3) move from right to left, you can perform the operation depicted in Figure \ref{e1-3m} and create a cone in that region.

\begin{figure}[htbp]
\begin{center}
\includegraphics[width=100mm]{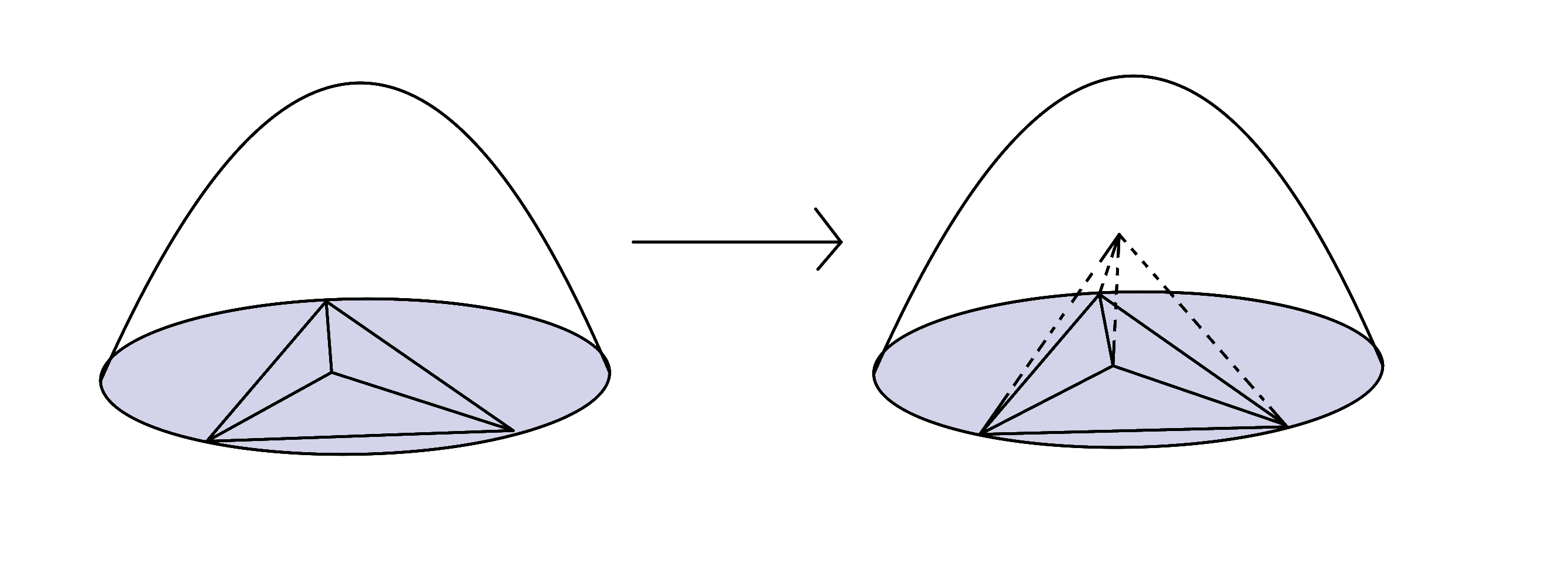}
\caption{Deformation near the boundary}
\label{e1-3m}
\end{center}
\end{figure}

Figure \ref{e1-3m} illustrates the concept of inserting a cone, shaped as shown, into the portion where the boundary is to be transformed. The resulting structure, with the inserted cone, remains a simplicial complex. By applying this operation once, you can transform the original simplicial complex $K'_a$ into a new one, denoted as $K'_b$. $|K'_a|=|K'_b|$, and the identity map $id: |K'_a|\to|K'_b|$ is a PL homeomorphism. Fundermore, since $id|_{\del K'_a}: \del K'_a\to \del K'_b$ is a simplicial map isomorphism, it follows from Thm \ref{th:tv} that $K'_a$ can be transformed into $K'_b$ using a three-dimensional Pachner move. Therefore, for the portion where you want to transform the two-dimensional Pachner move (1-3) move from right to left, you can consider that portion as a cone. By treating it as a cone, you can execute a three-dimensional (1-3) move and transform the boundary accordingly.
\\

{(b) About (2-2)move}

The reason why the portion transformed by the two-dimensional Pachner move (2-2) move can be replaced with a three-dimensional (2-2) move is as follows: For the portion where you want to perform the transformation using the two-dimensional Pachner move (2-2) move, you can apply the operation depicted in Figure \ref{e2-2m} and create a cone in that region.

\begin{figure}[htbp]
\begin{center}
\includegraphics[width=100mm]{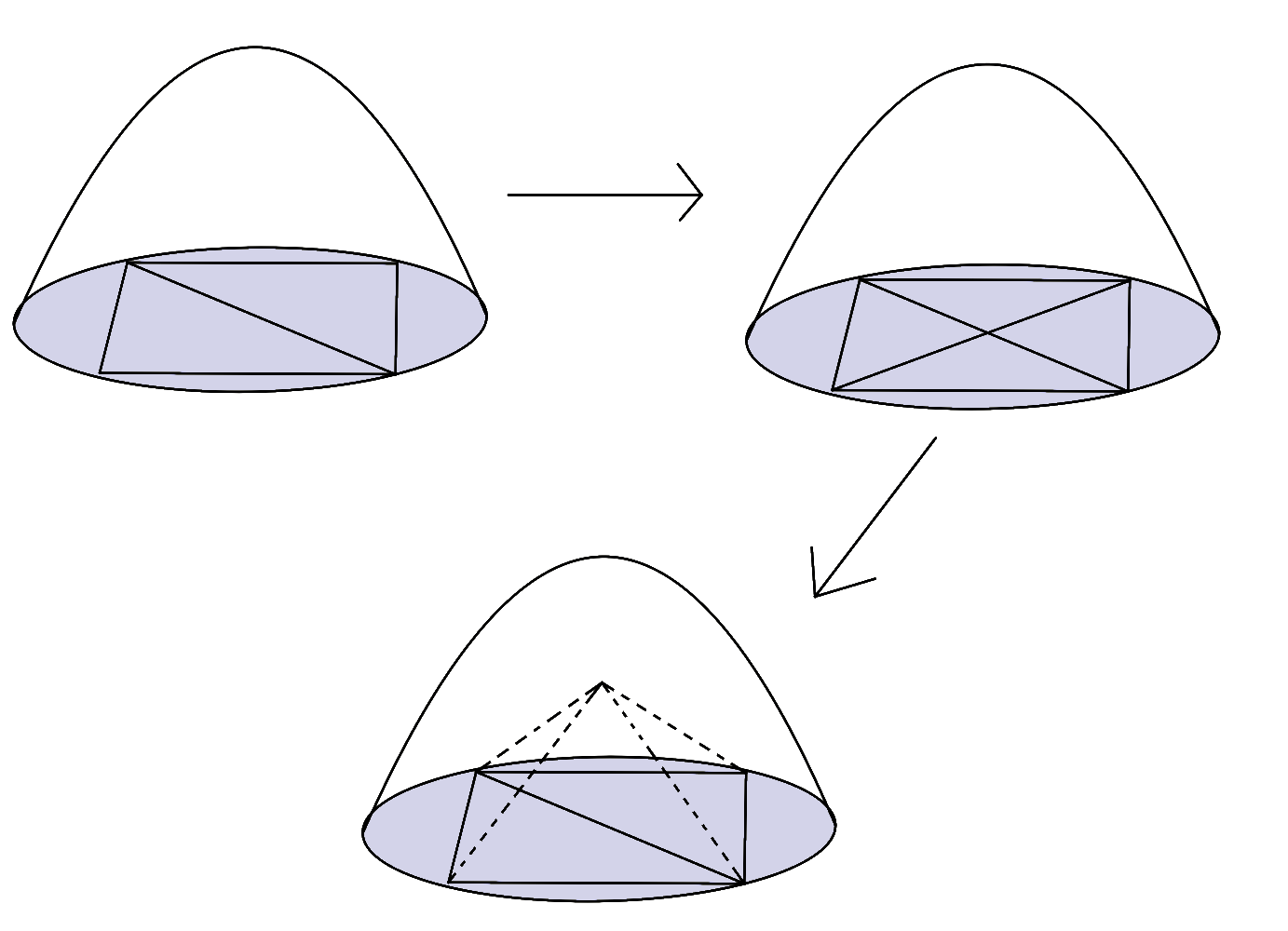}
\caption{Deformation near the boundary}
\label{e2-2m}
\end{center}
\end{figure}

First, replace the simplicial complex of the boundary region with a structure resembling a barycentric subdivision only for the region where the transformation is desired, while keeping the rest unchanged. In this case, according to Theorem \ref{th:btri}, it is possible to extend this simplicial subdivision of the boundary to consider the subdivision of the entire complex. Subsequently, Figure \ref{e2-2m} envisions inserting a cone into the region where the boundary is desired to be transformed, similar to Figure \ref{e1-3m}.

Even after inserting the cone, the structure remains a simplicial complex. By applying this sequence of operations once, you can transform the original simplicial complex $K'_a$ into a new one, denoted as $K'_b$.
The subsequent discussion is similar to the case (1). That is, $|K'_a|=|K'_b|$, $id: |K'_a|\to|K'_b|$ is a PL homeomorphism, and $id|_{\del K'_a}: \del K'_a\to \del K'_b$ is a simplicial map isomorphism. Therefore, $K'_a$ can be transformed into $K'_b$ using a three-dimensional Pachner move. Consequently, for the portion where you want to transform the two-dimensional Pachner move (2-2) move, you can consider that portion as a cone. By treating it as a cone, you can execute a three-dimensional (2-2) move and transform the boundary accordingly.
\\

From the above considerations, it follows that by using three-dimensional Pachner moves, three-dimensional (1-3) moves, and three-dimensional (2-2) moves, it is possible to make the triangular K coincide with L.
\end{proof}

\section{Example}

Here, we will show some of the result of using this invariant. In this section, we will give the result of 3 different basic geometris and 3 different complement of a knot. Knot we used is $4_1$ knot, torus knot (5,2). and $5_1$ knot.
\\

{\bf (1) Result of $D^3$}

\begin{figure}[htbp]
\begin{center}
\includegraphics[width=60mm]{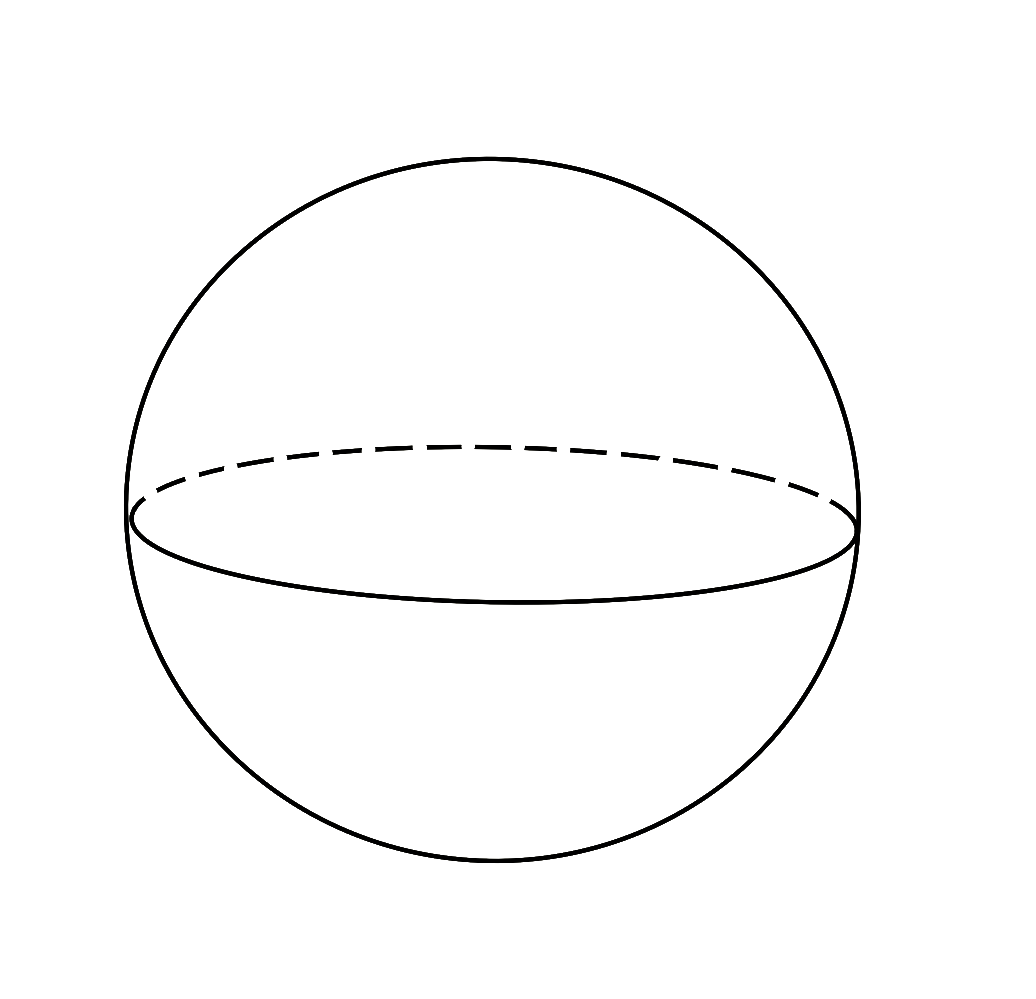}
\caption{Image of $D^3$}
\label{d3}
\end{center}
\end{figure}

\begin{figure}[htbp]
\begin{center}
\includegraphics[width=60mm]{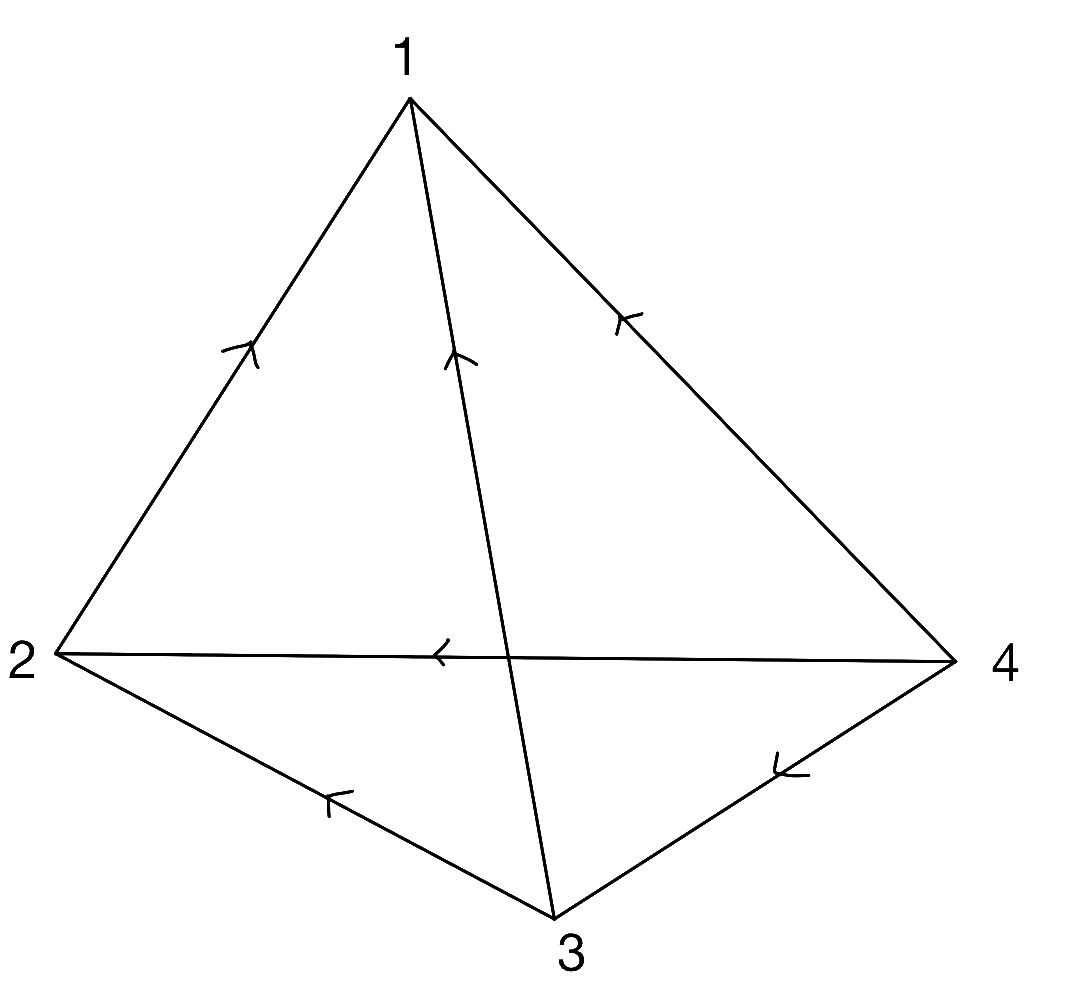}
\caption{Triangulation of $D^3$}
\label{d3tri}
\end{center}
\end{figure}

$D^3$ means next set:$D^3 \df \{(x,y,z)\in \mathbb{R}^3|\ x^2+y^2+z^2 \le1 \}$. This image will be like Figure \ref{d3} and triangular which I used will be like Figure \ref{d3tri}

As you can caluculate, color will be uniquely defined by label of 3 different 1-simplex and 3 different 2-simplex: $g_{12},g_{13},g_{14},h_{123},h_{134},h_{124}$. The invariant will be as follow:

\beq
\begin{split}
Z(D^3)&=|G|^{-4+6-4}|H|^{4-6+4-1}\\
&\quad\times{1\over |G|^6}\sum_{g_{12}\in G}\sum_{g_{13}\in G}\sum_{g_{14}\in G}\sum_{g_{23}\in G}\sum_{g_{24}\in G}\sum_{g_{34}\in G}\\
&\quad\quad\times{1\over |H|^4}\sum_{h_{123}\in H}\sum_{h_{124}\in H}\sum_{h_{134}\in H}\sum_{h_{234}\in H}\\
&\quad\quad\quad\times\delta_G(g_{123})\delta_G(g_{124})\delta_G(g_{134})\delta_G(g_{234})\delta_H(h_{1234})\\
&={1\over |G|^{2}}|H|{1\over |G|^6}{1\over |H|^4}\\
&\quad\times \sum_{g_{12}\in G}\sum_{g_{13}\in G}\sum_{g_{14}\in G}\sum_{h_{123}\in H}\sum_{h_{124}\in H}\sum_{h_{134}\in H}|G|^4|H|^1\\
&={1\over |G|^{2}}|H|{1\over |G|^6}{1\over |H|^4}|G|^3|H|^3|G|^4|H|^1\\
&={1\over |G|}|H|
\end{split}
\label{inv_d3}
\eeq
\\

{\bf (2) Result of $D^2\times S^1$}

\begin{figure}[htbp]
\begin{center}
\includegraphics[width=60mm]{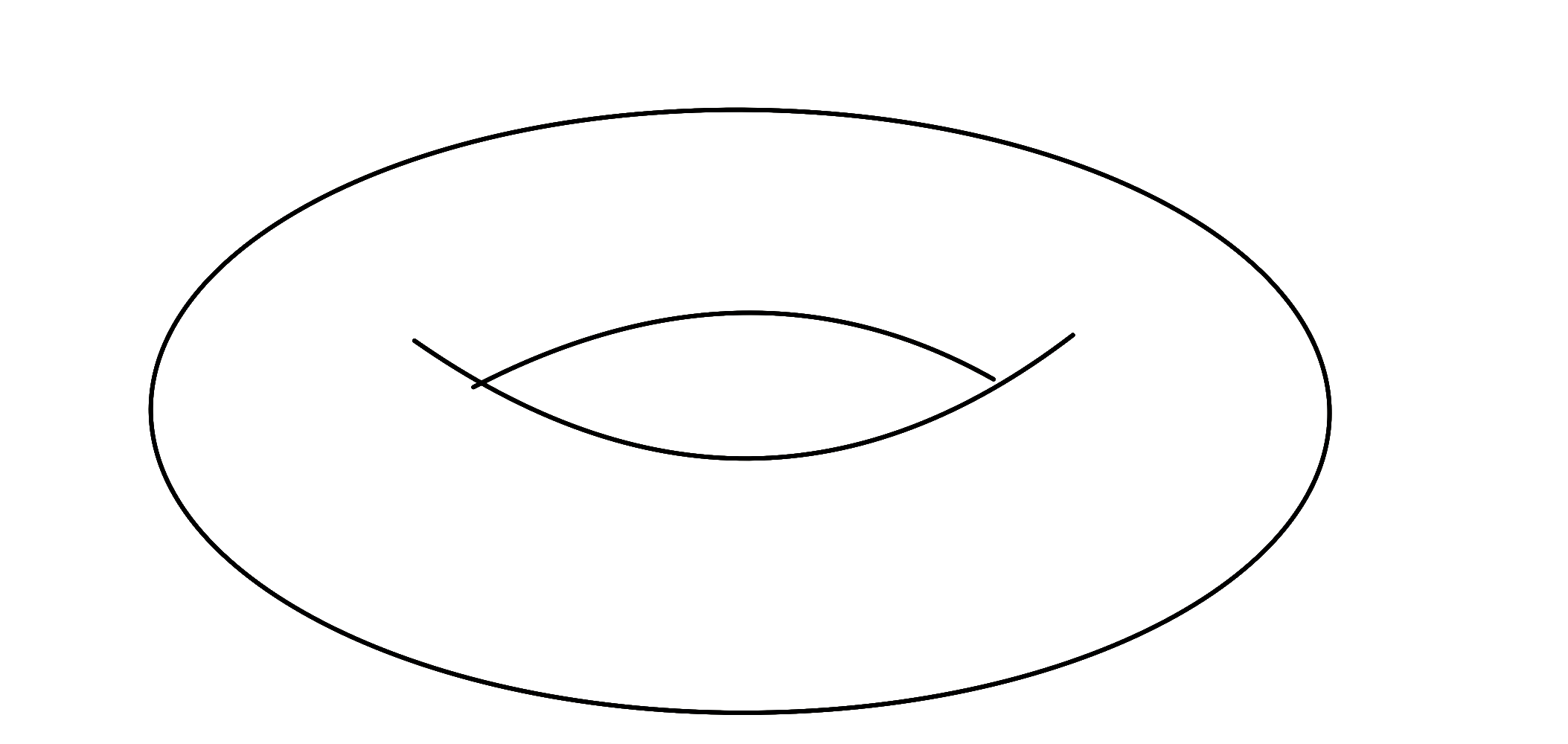}
\caption{Image of $D^2\times S^1$}
\label{d2s1}
\end{center}
\end{figure}

\begin{figure}[htbp]
\begin{center}
\includegraphics[width=60mm]{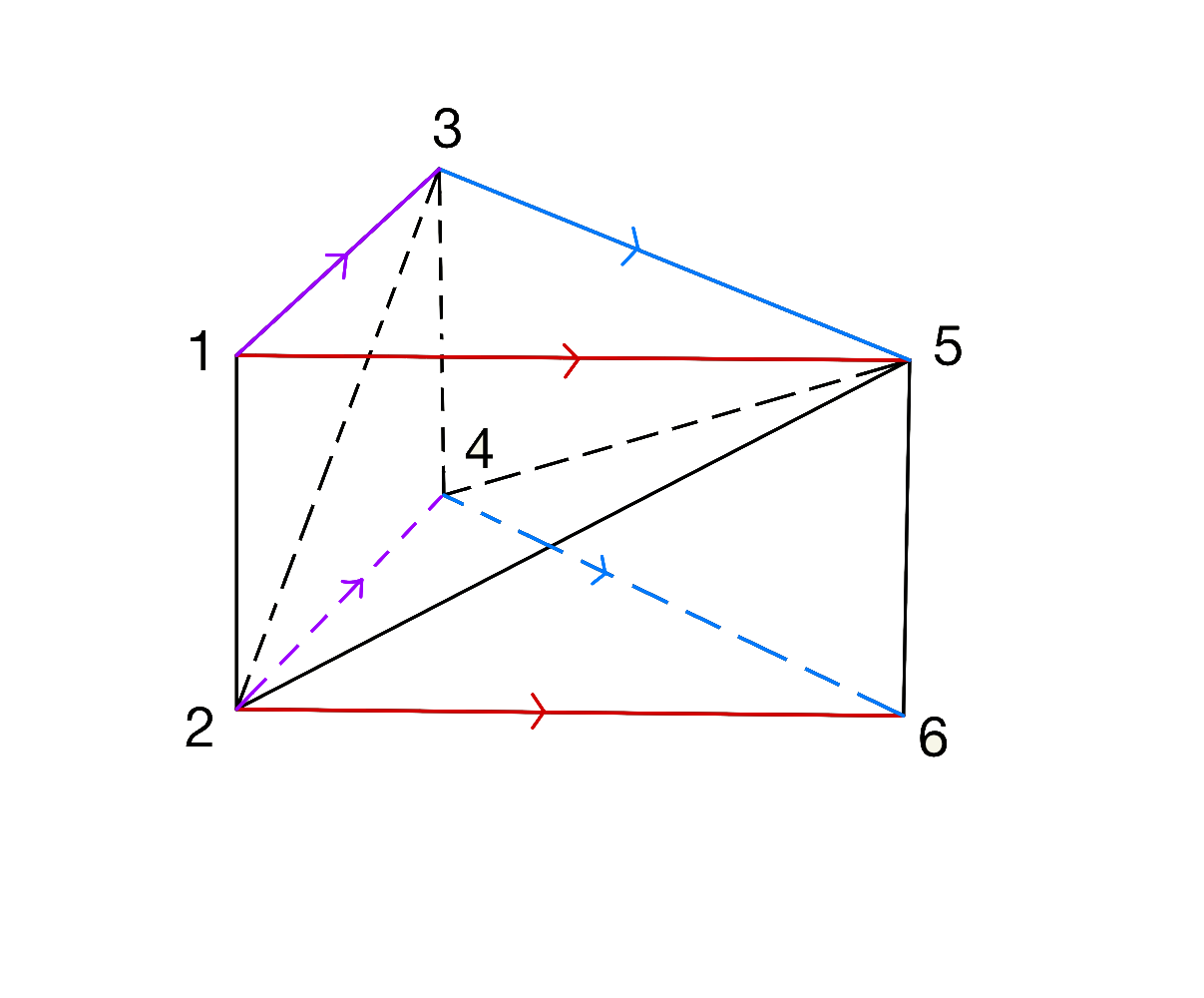}
\caption{Singular triangulation of $D^2\times S^1$}
\label{d2s1tri}
\end{center}
\end{figure}

The image of $D^2\times S^1$ will be like Figure \ref{d2s1} and Figure \ref{d2s1tri} will be an one of the example of singular triangulation of $D^2\times S^1$. You can check that this invariant can be caluculated by triangular with a local order. Because of this, we can used singular triangular instead of triangular.

The color of singular triangulation of $D^2\times S^1$ will be uniquely defined by value of $g_{12},g_{25},g_{46},h_{125},h_{256},h_{456},h_{345},h_{246}$. The invariant will be as follow:

\beq
\begin{split}
Z(D^2\times S^1)&=|G|^{-3+9-9}|H|^{3-9+9-3}\\
&\quad\times\left(\prod_{(jk)\in K_1}{1\over |G|}\sum_{g_{jk}\in G}\right)\left(\prod_{(jkl)\in K_2}{1\over |H|}\sum_{h_{jkl}\in H}\right)\\
&\quad\quad\times\left(\prod_{(jkl)\in K_2} \delta_G(g_{jkl})\right) \left(\prod_{(jklm)\in K_3} \delta_H(h_{jklm})\right)\\
&=|G|^{-3}{1\over |G|^9}{1\over |H|^9} \sum_{g_{12}\in G} \sum_{g_{25}\in G}\sum_{g_{46}\in G}\\
&\quad\times\sum_{h_{125}\in H}\sum_{h_{256}\in H}\sum_{h_{456}\in H}\sum_{h_{345}\in H}\sum_{h_{246}\in H}\sum_{h_{234}\in H}
|G|^9|H|^3\\
&=1
\end{split}
\label{inv_d2s1}
\eeq
\\

{\bf (3) Result of $S^2\times [0,1]$}

\begin{figure}[htbp]
\begin{center}
\includegraphics[width=60mm]{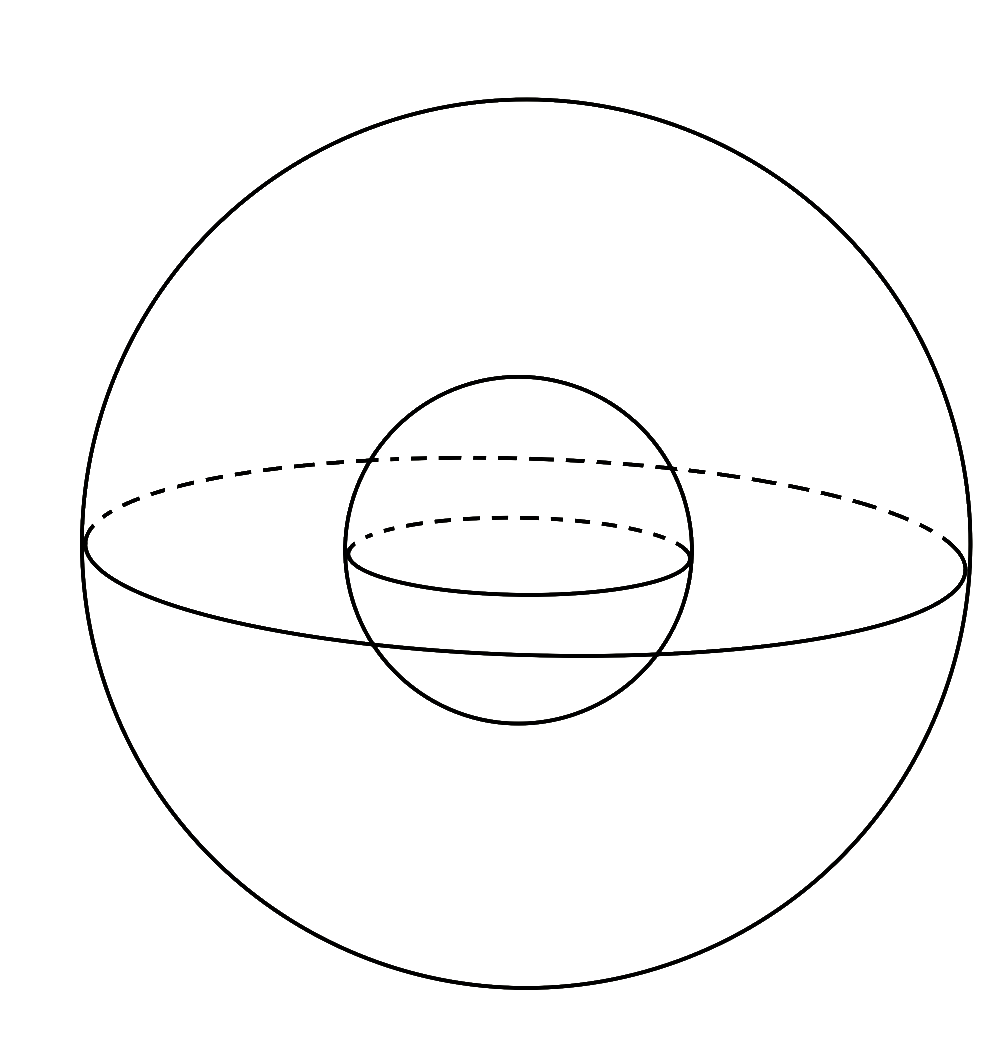}
\caption{Image of $S^2\times [0,1]$}
\label{s201}
\end{center}
\end{figure}

\begin{figure}[htbp]
\begin{center}
\includegraphics[width=60mm]{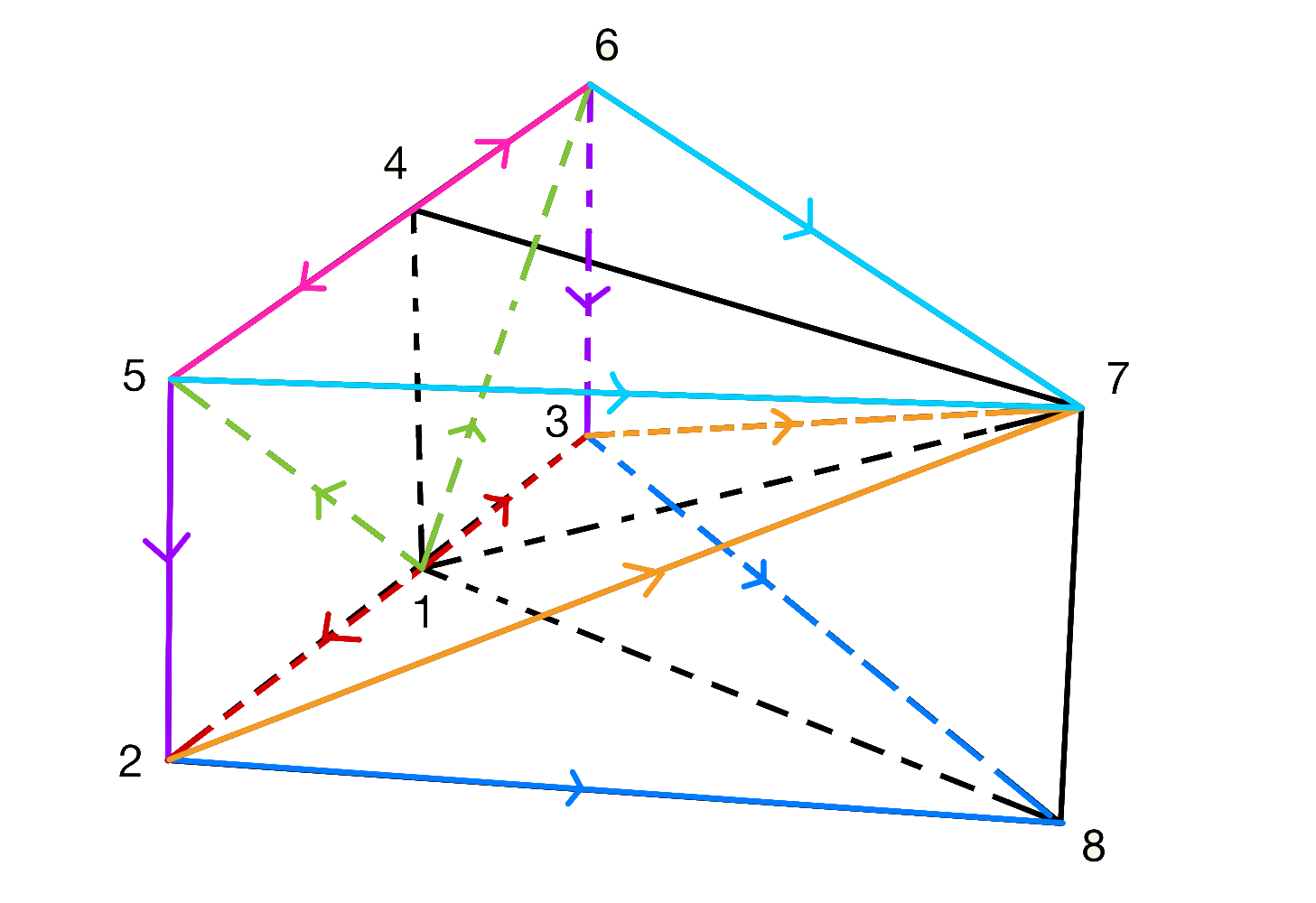}
\caption{Singular triangulation of $S^2\times [0,1]$}
\label{s201tri}
\end{center}
\end{figure}

The image of $S^2\times [0,1]$ will be like Figure \ref{s201} and Figure \ref{s201tri} will be an one of the example of singular triangulation of $S^2\times [0,1]$.

The color of singular triangulation of $S^2\times [0,1]$ will be uniquely defined by value of $g_{25},g_{57},g_{78},g_{18},g_{14},h_{257},h_{278},h_{128},h_{178},h_{125},h_{147},h_{145},h_{138}$. In this case, $h_{138}$ can not be taken for all $H$. Instead of it, $h_{138}$ have to satisfy following equation:$\del(h_{128})=\del(h_{138})$. The invariant will be as follow:
\beq
\begin{split}
Z(S^2\times [0,1])&=|G|^{-6+12-14}|H|^{6-12+14-6}\\
&\quad\times\left(\prod_{(jk)\in K_1}{1\over |G|}\sum_{g_{jk}\in G}\right)\left(\prod_{(jkl)\in K_2}{1\over |H|}\sum_{h_{jkl}\in H}\right)\\
&\quad\quad\times\left(\prod_{(jkl)\in K_2} \delta_G(g_{jkl})\right) \left(\prod_{(jklm)\in K_3} \delta_H(h_{jklm})\right)\\
&=|G|^{-8}|H|^2{1\over |G|^{12}}{1\over |H|^{14}} \sum_{g_{25}\in G} \sum_{g_{57}\in G}\sum_{g_{78}\in G}\sum_{g_{18}\in G}\sum_{g_{14}\in G}\sum_{h_{257}\in H}\sum_{h_{278}\in H}\sum_{h_{128}\in H}\\
&\quad\times\sum_{h_{178}\in H}\sum_{h_{125}\in H}\sum_{h_{147}\in H}\sum_{h_{145}\in H}\sum_{h_{138}\in H}
|G|^{13}\delta_G\left(\del(h^{-1}_{128}h_{138})\right)|H|^6\\
&={1\over |G|^8}|H|^2{1\over |G|^{12}}{1\over |H|^{14}}|G|^5|H|^6|G|^{13}|H|^6\sum_{h_{138}\in H}\sum_{h_{128}\in H}\delta_G(\del(h^{-1}_{128}h_{138}))\\
&={1\over|G|^2}\sum_{h_{138}\in H}\sum_{h_{128}\in H}\delta_G(\del(h^{-1}_{128}h_{138}))
\end{split}
\label{inv_s201}
\eeq

The difference between the value of invariant of $D^3$ and  $S^2\times [0,1]$ come from hall. When we fill $S^2\times [0,1]$, we will have a new 3-simplex in it so the $h_{138}$ in $S^2\times [0,1]$ have to be equal to
$h_{128}$. Then the equation (\ref{inv_s201}) will be coinside with equation (\ref{inv_d3}).

{\bf (4) Result of complement of $4_1$ knot}

\begin{figure}[htbp]
\begin{center}
\includegraphics[width=60mm]{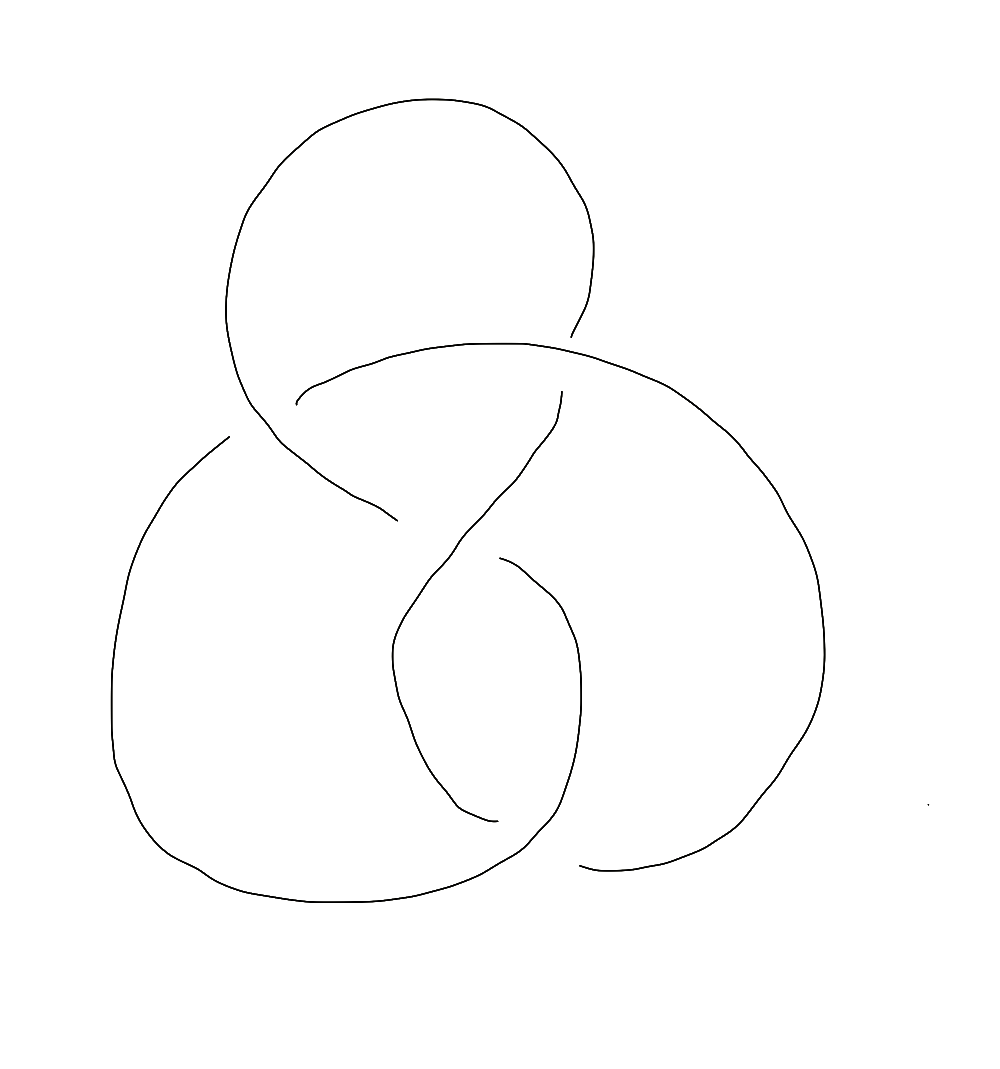}
\caption{Image of $4_1$ knot}
\label{k41}
\end{center}
\end{figure}

\begin{figure}[htbp]
\begin{center}
\includegraphics[width=60mm]{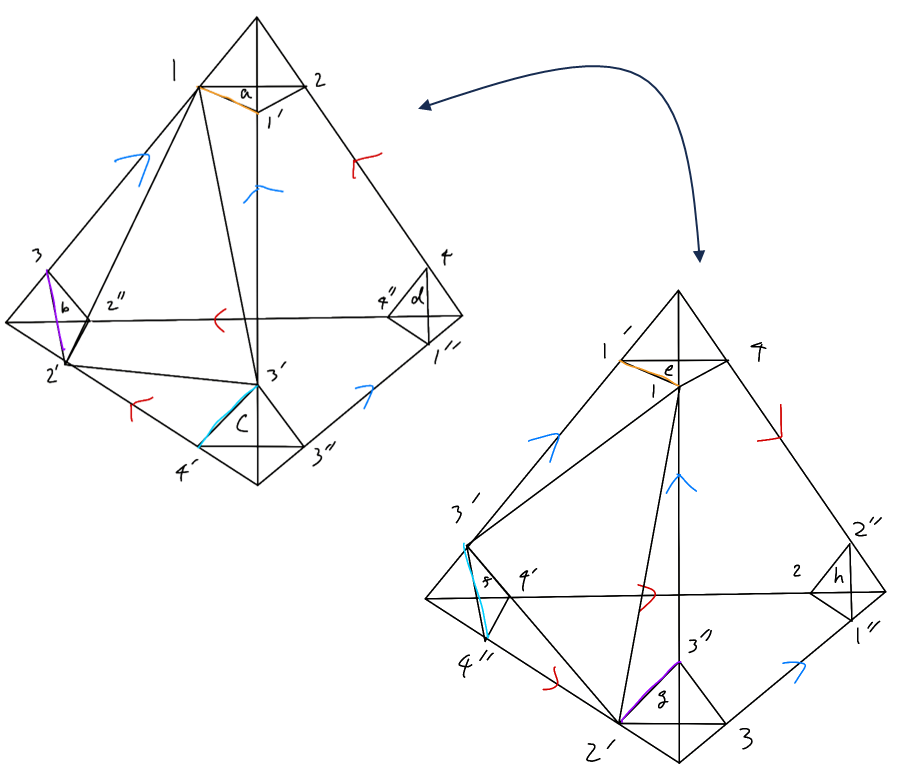}
\caption{Triangulation of complement of $4_1$ knot}
\label{k41tri}
\end{center}
\end{figure}

$4_1$ knot will be like Figure \ref{k41} and triangulation of complement of $4_1$ knot will be like Figure \ref{k41tri}. When illustrating the thought-out divisions, presenting only one side here to avoid clutter in the diagram. The same approach is taken for the other sides as well.

We will give the result first. The result will be as follow:
\beq
Z(S^3-4_1)={1\over |G|^2}{1\over |H|}\sum_{X,Y\in G, A\in H}\delta_G\left(XY^{-1}X^{-1}YX^{-1}Y^{-1}XYX^{-1}Y\del(A^{-1})\right)
\label{41zm}
\eeq

This equation holds for any crossed module $(H\stackrel{\del}{\to}G,\rhd)$ with arbitrary finite groups $H$ and $G$. However, fundamentally, it becomes apparent through calculations specifically when $H$ is a finite abelian group. In other words, the results are obtained by computing for all $h\in H$ with $\del(h)=e_G$.

\begin{figure}[htbp]
\begin{center}
\includegraphics[width=60mm]{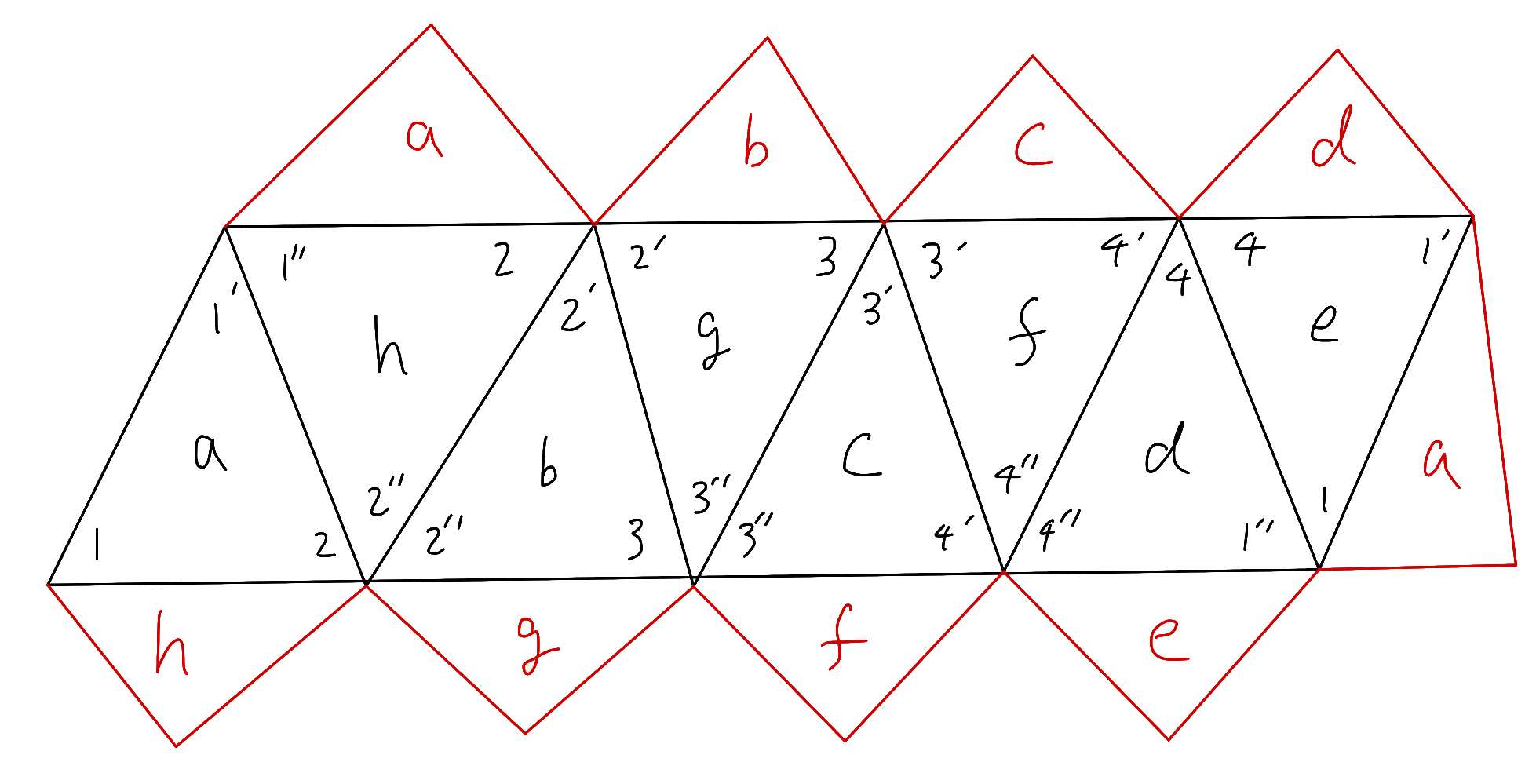}
\caption{Triangulation of boundary of complement of $4_1$ knot}
\label{41bou}
\end{center}
\end{figure}

To calucutate this first we have to think about surface of left tetrahedron.
The color of surface of left tetrahedron will be defined uniquely by $b,r,g_{11'},g_{3'4'},g_{2'2''},g_{3''4'},g_{1'2}$.$b,r$ is the color of 1-simplex with bule and red arrow. By gluing it to right  tetrahedron, the illustration depicting how it shares edges with other simplices is shown in Figure  \ref{41bou}. However, at this point, it is necessary to investigate the condition $g_{jkl}=e$ for the 2-simplex labeled with alphabets in the diagram. 

This leads to obtaining the conditions for the following equation \ref{41ebun}.
\begin{subequations} \label{41ebun}
\begin{empheq}[left = {\empheqlbrace \,}, right = {}]{align}
g_{1'2}&=g_{2'2''}g_{1'2}g_{11'} \label{41a} \\
bg^{-1}_{11'}b^{-1}g^{-1}_{3'4'}&=g^{-1}_{3''4'}g_{3'4'}bg^{-1}_{11'}b^{-1}g^{-1}_{3'4'}g^{-1}_{2'2''} \label{41b} \\
g_{3'4'} &= rg_{2'2''}r^{-1}g_{3'4'}bg^{-1}_{1'2}r^{-1}g_{3''4'}^{-1}  \label{41c} \\
rg_{1'2}b^{-1}g_{3'4'}^{-1}g_{3''4}b&=rg_{2'2''}r^{-1}g_{3''4'}bg_{11'} \label{41d} 
\end{empheq}
\end{subequations}
Indeed, upon examining this equation, it becomes evident that when equations (\ref{41a})-(\ref{41c}) hold true, equation (\ref{41d}) also holds. Therefore, solving equations (\ref{41a})-(\ref{41c}) is sufficient. Solving these equations yields equation (\ref{41fin}). Moreover, when equation (\ref{41fin}) holds, unique values for $g_{2'2''}$ and $g_{1'2}$ are determined such that equations (\ref{41a})-(\ref{41c}) are satisfied.
\beq
\begin{split}
e &=g_{3'4'}^{-1}g_{3''4}bg_{11'}^{-1}b^{-1}g_{3''4}^{-1}g_{3'4'}bg_{11'}b^{-1}g_{3''4}^{-1}g_{3'4'}bg_{11'}^{-1}b^{-1}g_{3'4'}^{-1}g_{3''4'}bg_{11'}b^{-1}g_{3''4}^{-1}g_{3'4'}bg_{11'}b^{-1}\\
& =XY^{-1}X^{-1}YX^{-1}Y^{-1}XYX^{-1}Y
\end{split}
\label{41fin}
\eeq
However, by transforming variables as $X\df g_{3'4'}^{-1}g_{3''4}$ and $Y\df bg_{11'}b^{-1}$, we arrive at an expression akin to equation (\ref{41zm}). Do a same thing with in general Crossed module we get equation (\ref{41zm}).

{\bf (5) Result of complement of torus knot (5,2)}

\begin{figure}[htbp]
\begin{center}
\includegraphics[width=60mm]{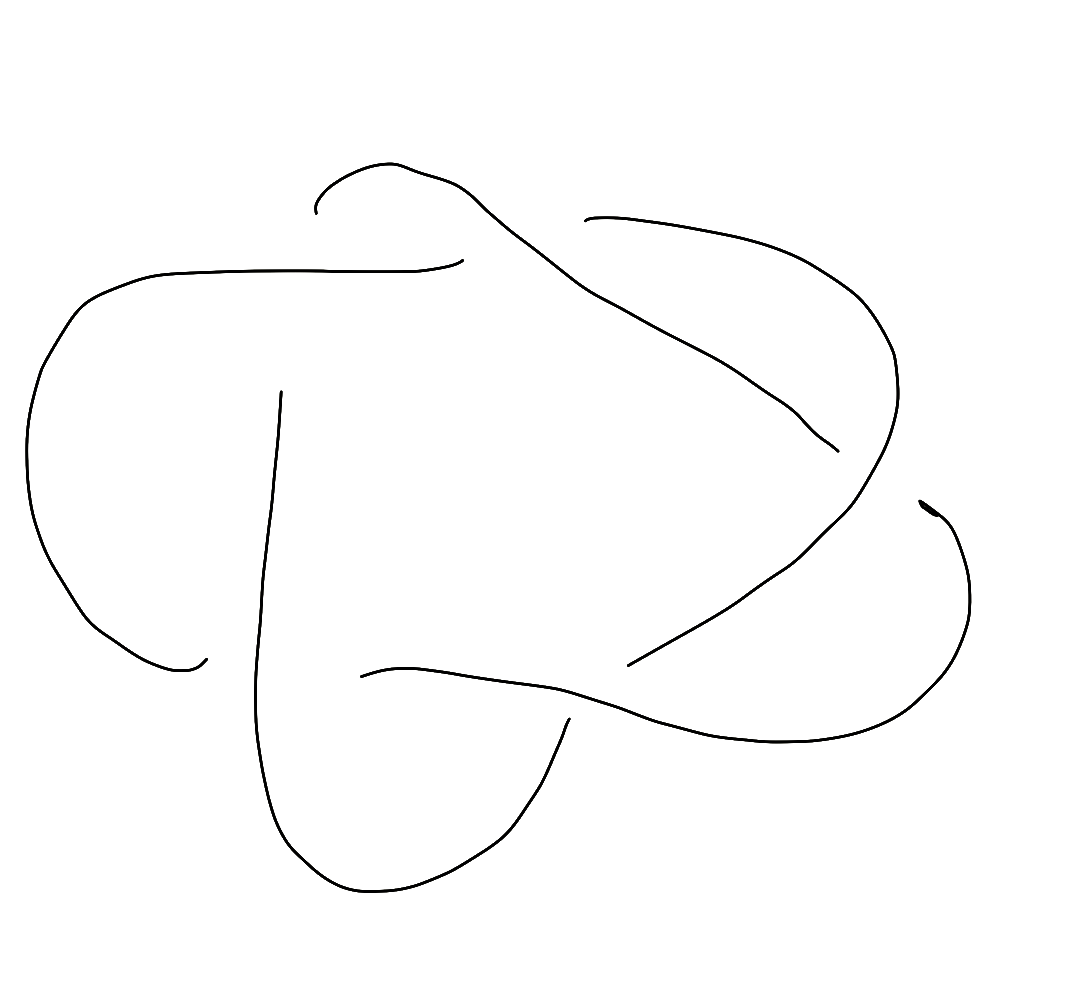}
\caption{Image of torus knot (5,2)}
\label{t52}
\end{center}
\end{figure}

\begin{figure}[htbp]
\begin{center}
\includegraphics[width=60mm]{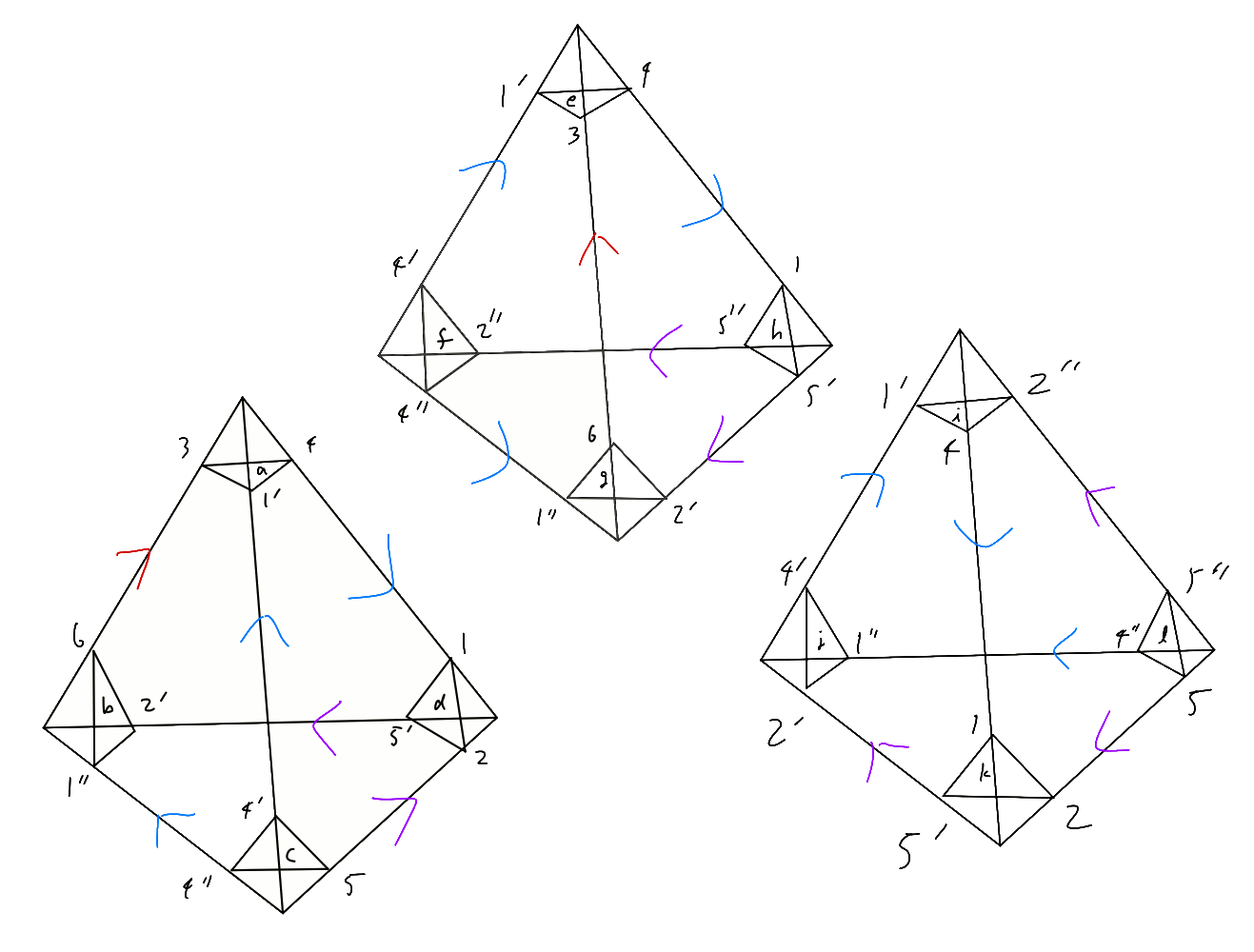}
\caption{Triangulation of complement of torus knot (5,2)}
\label{t52tri}
\end{center}
\end{figure}

This can be caluculate almost same as complement of $4_1$ knot. We just write the result.
Torus knot (5,2) is like Figure \ref{t52} and triangulation of complement of torus knot (5,2) is look like Figure \ref{t52tri}.

The result will be as follow:
\beq
Z(S^3-T(5,2))={1\over |G|^2}{1\over |H|}\sum_{X,Y\in G, A\in H}\delta_G\left(YX^4YX^{-1}\del(A^{-1})\right)
\label{t52zm}
\eeq

{\bf (6) Result of complement of $5_2$ knot}

\begin{figure}[htbp]
\begin{center}
\includegraphics[width=60mm]{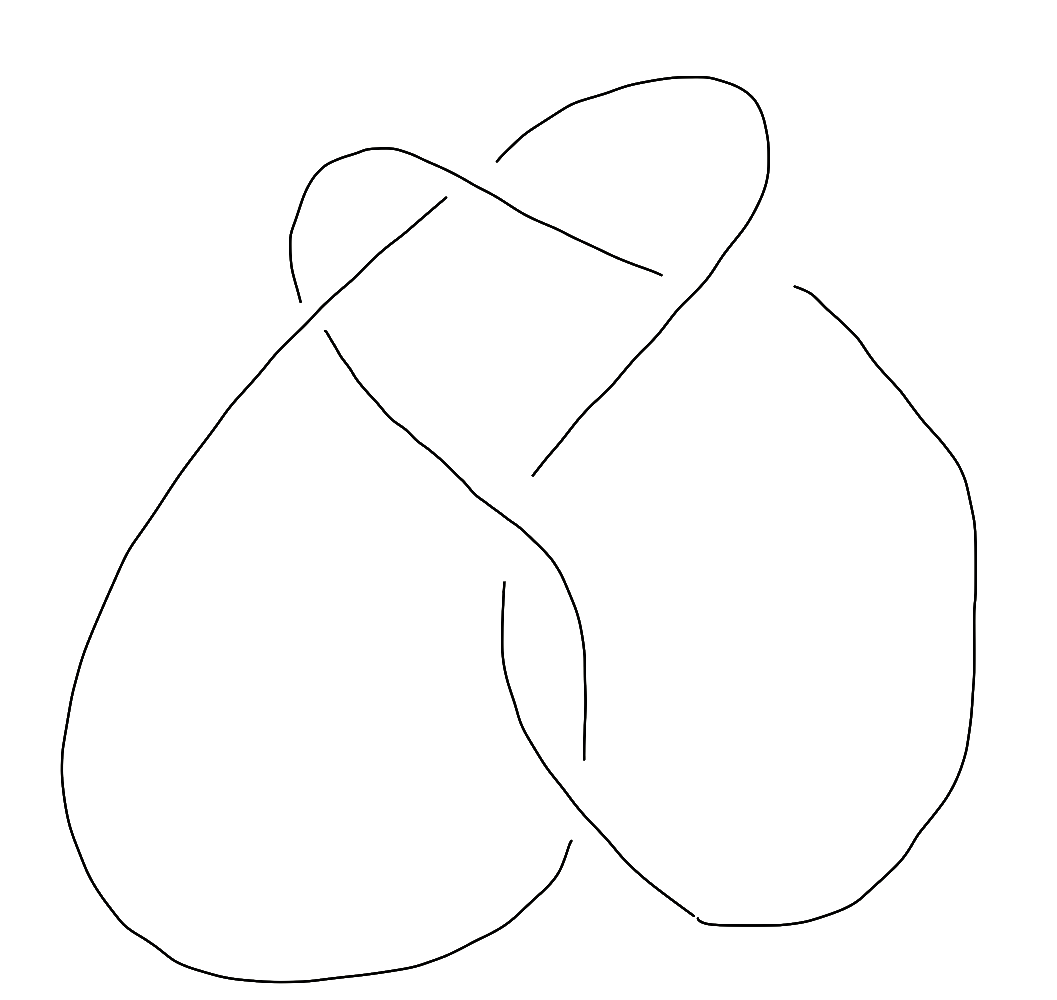}
\caption{Image of $5_2$ knot }
\label{k52}
\end{center}
\end{figure}

\begin{figure}[htbp]
\begin{center}
\includegraphics[width=60mm]{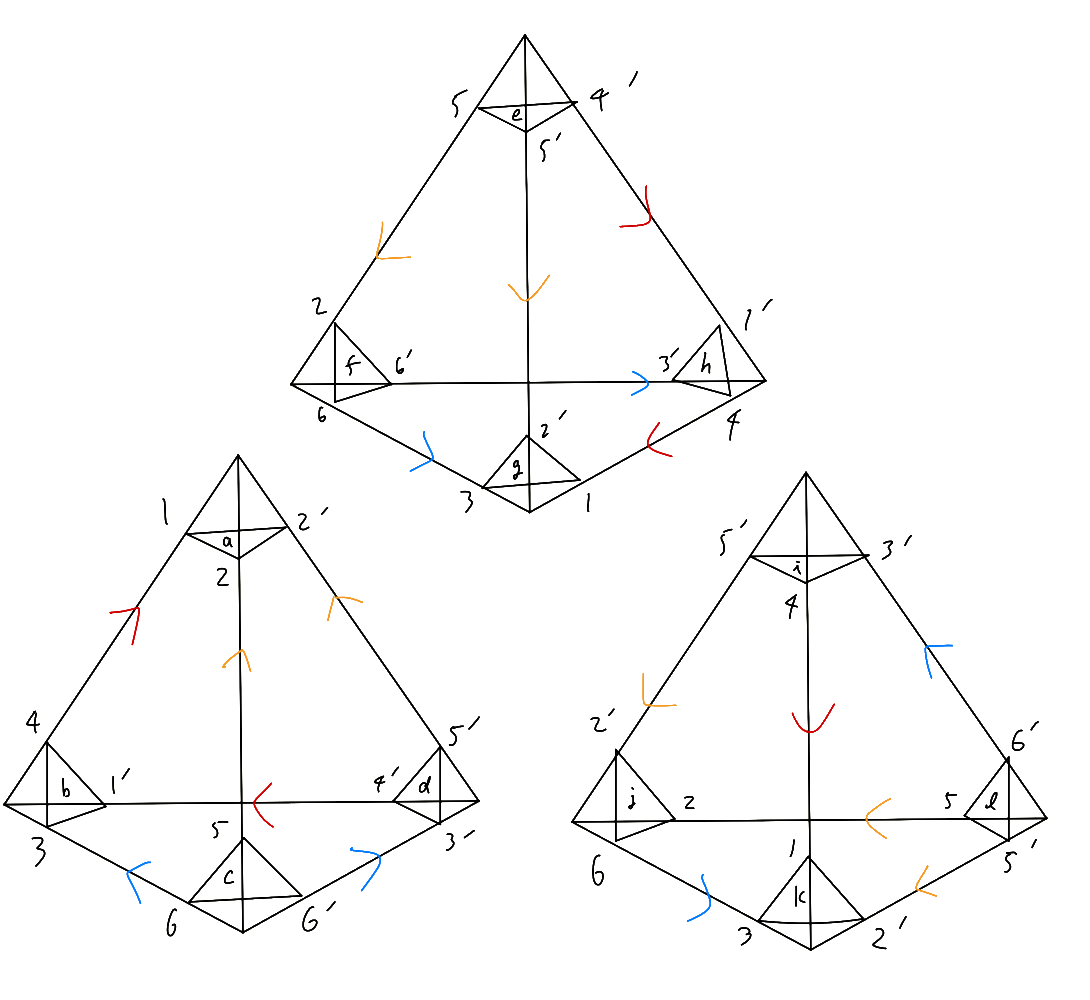}
\caption{Triangulation of complement of $5_2$ knot }
\label{k52tri}
\end{center}
\end{figure}

Complement of $5_2$ knot can be caluculated as same as other complement spaces of knot,so we just write the result.
$5_2$ knot is like Figure \ref{k52} and triangulation of complement of torus knot (5,2) is look like Figure \ref{k52tri}.

The result will be as follow:
\beq
Z(S^3-5_2)={1\over |G|^2}{1\over |H|}\sum_{X,Y\in G, A\in H}\delta_G\left(XY^{-1}X^{-1}YXY^{-1}X^{-1}YX^{-1}Y^{-1}XYX^{-1}Y^{-1}\del(A^{-1})\right)
\label{52kzm}
\eeq

We have examined three examples of the complement of knots, and it is evident that the content of $\delta_G$ is closely related to certain aspects. Specifically, it corresponds to the relations in the fundamental group of the complement of knots. In fact, for the $4_1$ knot and the Torus(5,2) knot, it has been verified that, apart from $\del(A^{-1})$, the content of $\delta_G$ is equivalent to the relations in the representation of the fundamental group. Therefore, even in the case of a manifold with boundary, when $H$ is trivial, it can be anticipated that the number of representations from the fundamental group of the manifold with boundary to $G$ is a constant multiplied by $1\slash|G|$.

\appendix

\section{3-dimensional (1-3)move and 3-dimensional (2-2)move}

In this chapter, we show that 3-dimensional (1-3)move and 3-dimensional (2-2)move preseve the Formula \ref{3mfd_w_b thm}.

\noindent
{\bf (a) 3-dimensional (1-3) move in Figure \ref{1-3mv} }
\\

First, let's examine the invariance of the equation (\ref{ZM}) for the 3-dimensional (1-3) move shown in Figure \ref{1-3mv}. Since it is known that this invariance is independent of the choice of the total order on the 0-simplices, we will perform the calculation using the total order depicted in Figure \ref{1-3mv}. We denote the diagram on the left side of Figure \ref{1-3mv} as l.h.s and the one on the right side as r.h.s. First, for the l.h.s, the equation (\ref{ZM}) can be expressed as follows:

\beq
\begin{split}
Z(l.h.s)&=|G|^{-4+6-4}|H|^{4-6+4-1}\left({1\over |H|}\sum_{h_{345}\in H}\right)\delta_G(g_{345})\delta_H(h_{2345})Z_{remainder}\\
&=|G|^{-2}|H|
\left({1\over |H|}
	{ \sum_{h_{345}\in H}}
	\right)\\
& \quad\quad\times 
\delta_G\Big(
\del(h_{345})g_{45}g_{34}g_{35} ^{-1}
\Big)\\
&\quad\quad\quad\times \delta_H\Big(
h_{245}(g_{45}\rhd h_{234})h_{345}^{-1}h_{235}^{-1}\Big)
Z_{remainder}\\
\end{split}
\label{13lhs}
\eeq

On the other hand, when using the equation (\ref{ZM}) for the r.h.s, it can be expressed as follows:
\beq
\begin{split}
Z(r.h.s)&=|G|^{-5+10-9}|H|^{5-10+9-3}
\left(
	 \prod_{(jk) \in M_1}{1\over |G|} \sum_{g_{jk}\in G}
\right)
\left(
	 \prod_{(jkl) \in M_2}{1\over |H|} \sum_{h_{jkl}\in H}
\right)\\
& \quad\times
\delta_G(g_{123})\delta_G(g_{124})\delta_G(g_{125})\delta_G(g_{134})\delta_G(g_{135})\delta_G(g_{145})\\
& \quad\quad\times\delta_H(h_{1234})\delta_H(h_{1235})\delta_H(h_{1245})
Z_{remainder}
\end{split}
\label{13rhs_1}
\eeq
Here, $Z_{\text{remainder}}$ refers to the part of the sum and the common $g_{jkl}$, $h_{jklm}$ terms shared between l.h.s and r.h.s. The sets $M_1$ and $M_2$ in r.h.s are defined as follows:
\[
\begin{split}
M_1&=\{(12), (13), (14), (15) \},\\
M_2&=\{(123), (124), (125), (134), (135), (145)\}, \\
\end{split}
\]
Let's consider the transformation of the r.h.s. First, for the 2-simplex $|123|$ in r.h.s, since $\delta_G(g_{123})$, when the coloring is correct, is equal to $e$, we can use this to transform it as follows:

\beq
\begin{split}
g_{123}&\df \del(h_{123})g_{23}g_{12}g_{13}^{-1} = e\\
&\iff g_{12}= g_{23}^{-1}\del(h_{123}^{-1})g_{13}
\end{split}
\label{13g12}
\eeq
By restricting $g_{12}$ to satisfy the equation (\ref{13g12}), $\delta_G(g_{123})$ will always be equal to $|G|$. Furthermore, from the equation (\ref{13g12}) and the coloring restriction on the 3-simplex $|1234|$, it follows that $h_{1234}=e$, leading to $g_{124}$ always being equal to $e$. Substituting the result of equation (\ref{13g12}) into $g_{12}$ yields the following:
\beq
\begin{split}
g_{124}&\df \del(h_{124})g_{24}g_{12}g_{14}^{-1} \\
&= \del(h_{124})g_{24}g_{23}^{-1}\del(h_{123}^{-1})g_{13}g_{14}^{-1}
\end{split}
\label{13g124_1}
\eeq
Furthermore, from the coloring restriction on the 3-simplex $|1234|$, the following holds due to $h_{1234}=e$:
\beq
\begin{split}
h_{1234}&\df h_{134}(g_{34}\rhd h_{123})h_{234}^{-1}h_{124}^{-1}=e \\
&\iff h_{123}^{-1} = g_{34}^{-1} \rhd (h_{234}^{-1}h_{124}^{-1}h_{134})
\end{split}
\label{13h1234_123}
\eeq
Substituting the result of equation (\ref{13h1234_123}) into $h_{123}^{-1}$ in equation (\ref{13g124_1}), we get the following:
\beq
\begin{split}
g_{124}&= \del(h_{124})g_{24}g_{23}^{-1}\del(h_{123}^{-1})g_{13}g_{14}^{-1}\\
&=\del(h_{124})g_{24}g_{23}^{-1}\del\Big(g_{34}^{-1} \rhd (h_{234}^{-1}h_{124}^{-1}h_{134})\Big)g_{13}g_{14}^{-1}\\
&=\del(h_{124})g_{24}g_{23}^{-1}g_{34}^{-1}\del(h_{234}^{-1})\del(h_{124}^{-1})\del(h_{134})g_{34}g_{13}g_{14}^{-1}
\end{split}
\label{13g124_2}
\eeq
Here, from the coloring restriction $\delta_G(g_{234})$ on the 2-simplex $|234|$, when the coloring is correct, $g_{234}=e$. Similarly, from the coloring restriction $\delta_G(g_{134})$ on the 2-simplex $|134|$, when the coloring is correct, $g_{134}=e$. These two conditions imply the following equation:
\beq
\begin{split}
g_{234}&\df \del(h_{234})g_{34}g_{23}g_{24}^{-1} =e \\
&\iff \redunderline{g_{24}g_{23}^{-1}g_{34}^{-1}\del(h_{234}^{-1})=e} \\
g_{134}&\df \del(h_{134})g_{34}g_{13}g_{14}^{-1} =e \\
&\iff \blueunderline{\del(h_{134})g_{34}g_{13}=g_{14}}
\end{split}
\label{13_234134}
\eeq
Using the equation (\ref{13_234134}), we can transform the equation (\ref{13g124_2}) as follows:
\beq
\begin{split}
g_{124}
&=\del(h_{124})\redunderline{g_{24}g_{23}^{-1}g_{34}^{-1}\del(h_{234}^{-1})}\del(h_{124}^{-1})\blueunderline{\del(h_{134})g_{34}g_{13}}g_{14}^{-1}\\
&=\del(h_{124})\redunderline{e}\del(h_{124}^{-1})\blueunderline{g_{14}}g_{14}^{-1}\\
&=e
\end{split}
\label{13g124_f}
\eeq
Therefore, by introducing $g_{12}$ as in equation (\ref{13g12}) and demonstrating that $g_{124}$ is always equal to $e$ when the other colors are correct, it has been shown.

Similarly, by introducing $g_{12}$ as in equation (\ref{13g12}) and demonstrating that $g_{125}$ is always equal to $e$ when the other colors are correct, we proceed. Eliminating $g_{12}$ in terms of $g_{125}$ using equation (\ref{13g12}), we get the following:
\beq
\begin{split}
g_{125}&\df \del(h_{125})g_{25}g_{12}g_{15}^{-1} \\
&= \del(h_{125})g_{25}g_{23}^{-1}\del(h_{123}^{-1})g_{13}g_{15}^{-1}
\end{split}
\label{13g125_1}
\eeq
Here, there is a coloring restriction $\delta_H(h_{1235})$ on the 3-simplex $|1235|$, and when the coloring is correct, $h_{1235} = e$. This implies the following:
\beq
\begin{split}
h_{1235}&\df h_{135}(g_{35}\rhd h_{123})h_{235}^{-1}h_{125}^{-1}=e \\
&\iff h_{123}^{-1} = g_{35}^{-1} \rhd (h_{235}^{-1}h_{125}^{-1}h_{135})
\end{split}
\label{13h1235_123}
\eeq
Substituting the result of equation (\ref{13h1235_123}) into $h_{123}^{-1}$ in equation (\ref{13g125_1}), we get the following:
\beq
\begin{split}
g_{125}
&= \del(h_{125})g_{25}g_{23}^{-1}\del(h_{123}^{-1})g_{13}g_{15}^{-1}\\
&=\del(h_{125})g_{25}g_{23}^{-1}\del\Big(g_{35}^{-1} \rhd (h_{235}^{-1}h_{125}^{-1}h_{135})\Big)g_{13}g_{15}^{-1}\\
&=\del(h_{125})g_{25}g_{23}^{-1}g_{35}^{-1}\del (h_{235}^{-1})\del (h_{125}^{-1})\del (h_{135})g_{35}g_{13}g_{15}^{-1}
\end{split}
\label{13g125_2}
\eeq
Here, from the coloring restriction $\delta_G(g_{235})$ on the 2-simplex $|235|$, when the coloring is correct, $g_{235}=e$. Similarly, from the coloring restriction $\delta_G(g_{135})$ on the 2-simplex $|135|$, when the coloring is correct, $g_{135}=e$. These two conditions imply the following equation:
\beq
\begin{split}
g_{235}&\df \del(h_{235})g_{35}g_{23}g_{25}^{-1} =e \\
&\iff \redunderline{g_{25}g_{23}^{-1}g_{35}^{-1}\del(h_{235}^{-1})=e} \\
g_{135}&\df \del(h_{135})g_{35}g_{13}g_{15}^{-1} =e \\
&\iff \blueunderline{\del(h_{135})g_{35}g_{13}=g_{15}}
\end{split}
\label{13_235135}
\eeq
Using the equation (\ref{13_235135}), we can transform the equation (\ref{13g125_2}) as follows:
\beq
\begin{split}
g_{125}
&=\del(h_{125})\redunderline{g_{25}g_{23}^{-1}g_{35}^{-1}\del(h_{235}^{-1})}\del(h_{125}^{-1})\blueunderline{\del(h_{135})g_{35}g_{13}}g_{15}^{-1}\\
&=\del(h_{125})\redunderline{e}\del(h_{125}^{-1})\blueunderline{g_{15}}g_{15}^{-1}\\
&=e
\end{split}
\label{13g125_f}
\eeq

Therefore, by introducing $g_{12}$ as in equation (\ref{13g12}) and demonstrating that $g_{125}$ is always equal to $e$ when the other colors are correct, it has been shown.

Based on what has been established so far, organizing the equation (\ref{13rhs_1}) yields the following:

\beq
\begin{split}
Z(r.h.s)&=|G|^{-4}|H|
\blueunderline{\left(
	 {1\over |G|} \sum_{g_{12}\in G}
\right)}
\left(
	 \prod_{(jk) \in M_1'}{1\over |G|} \sum_{g_{jk}\in G}
\right)
\left(
	 \prod_{(jkl) \in M_2}{1\over |H|} \sum_{h_{jkl}\in H}
\right)\\
& \quad\quad\times\redunderline{\delta_G(g_{123})\delta_G(g_{124})\delta_G(g_{125})}\delta_G(g_{134})\delta_G(g_{135})\delta_G(g_{145})\\
&\quad\quad\quad\times \delta_H(h_{1234})\delta_H(h_{1235})\delta_H(h_{1245})
Z_{remainder}\\
&=|G|^{-4}|H|
\blueunderline{{1\over |G|}}
\left(
	 \prod_{(jk) \in M_1'}{1\over |G|} \sum_{g_{jk}\in G}
\right)
\left(
	 \prod_{(jkl) \in M_2}{1\over |H|} \sum_{h_{jkl}\in H}
\right)\\
& \quad\quad\times\redunderline{ |G|^3} \delta_G(g_{134})\delta_G(g_{135})\delta_G(g_{145})\\
&\quad\quad\quad\times \delta_H(h_{1234})\delta_H(h_{1235})\delta_H(h_{1245})
Z_{remainder}\\
\end{split}
\label{13rhs_2}
\eeq
However, setting $M_1' ={(13),(14),(15)}$.

Here, the reason for being able to simplify as indicated by the blue underlines is that, due to equation (\ref{13g12}), once $g_{13}$, $g_{14}$, and $h_{123}$ are determined, $\delta_G(g_{123})$ uniquely determines $g_{12}$. The part underlined in red is always $g_{123}=g_{124}=g_{125}=e$ based on the results obtained so far. Comparing equation (\ref{13rhs_2}) with equation (\ref{13lhs}), apart from the summation part, the terms $\delta_G(g_{345})$ on the l.h.s and $\delta_G(g_{134})\delta_G(g_{135})\delta_G(g_{145})$ on the r.h.s are roughly equal. To establish the equality, it suffices to show that the terms $\delta_H(h_{2345})$ on the l.h.s and $\delta_H(h_{1234})\delta_H(h_{1235})\delta_H(h_{1245})$ on the r.h.s are also equal.

Therefore, to demonstrate the equality of $\delta_G(g_{345})$ on the l.h.s and $\delta_G(g_{134})\delta_G(g_{135})\delta_G(g_{145})$ on the r.h.s, consider the transformation of $\delta_G(g_{134})\delta_G(g_{135})\delta_G(g_{145})$ on the r.h.s.

For the restriction $g_{134}, g_{135}$ on the 2-simplex $|134|, |135|$ in r.h.s, when the colors are correct and $g_{134}=g_{135}=e$, it can be written as follows:
\beq
\begin{split}
g_{134} &\df \del(h_{134})g_{34}g_{13}g_{14}^{-1} = e\\
&\iff g_{14} = \del(h_{134})g_{34}g_{13} \\
g_{135} &\df \del(h_{135})g_{35}g_{13}g_{15}^{-1} = e\\
&\iff g_{15}^{-1} = g_{13}^{-1}g_{35} ^{-1}\del(h_{135}^{-1})\\
\end{split}
\label{13134135}
\eeq
Therefore, assuming that $g_{14}$ and $g_{15}$ always vary in such a way as to satisfy the equation (\ref{13134135}), substituting the result of equation (\ref{13134135}) into $g_{14}$ and $g_{15}$ in $\delta_G(g_{145})$ gives the following:
\beq
\begin{split}
\delta_G(g_{145}) &\df \delta_G\Big(\del(h_{145})g_{45}g_{14}g_{15}^{-1}\Big)\\
&=   \delta_G
\Big(
\del(h_{145})\blueunderline{g_{45}\del(h_{134})}g_{34}g_{13}g_{13}^{-1}g_{35} ^{-1}\redunderline{\del(h_{135}^{-1}})
\Big)  \\
&=\delta_G
\Big(
\redunderline{\del(h_{135}^{-1})}\del(h_{145})\blueunderline{\big(g_{45}\rhd\del(h_{134})\big)g_{45}}g_{34}g_{35} ^{-1}
\Big)\\
&=\delta_G
\Big(
\del(h_{345}')g_{45}g_{34}g_{35} ^{-1}
\Big)
\end{split}
\label{13r145}
\eeq
Here, we set $h_{345}'\df h_{135}^{-1}h_{145}\big(g_{45}\rhd h_{134}\big)$.

Next, consider the transformation of $\delta_H(h_{1234})\delta_H(h_{1235})\delta_H(h_{1245})$ on the r.h.s. Here, for the restriction $h_{1234}, h_{1235}$ on the 3-simplex $|1234|, |1235|$ in the r.h.s, when the colors are correct and $h_{1234}=h_{1235}=e$, it can be written as follows:
\beq
\begin{split}
h_{1234} &\df h_{134}(g_{34}\rhd h_{123})h_{234}^{-1}h_{124}^{-1} = e\\
&\iff h_{124}=h_{134}(g_{34}\rhd h_{123})h_{234}^{-1} \\
h_{1235} &\df h_{135}(g_{35}\rhd h_{123})h_{235}^{-1}h_{125}^{-1} = e\\
&\iff h_{125}^{-1} = h_{235}(g_{35}\rhd h_{123}^{-1})h_{135}^{-1}\\
\end{split}
\label{1312341235}
\eeq
Therefore, assuming that $h_{124}$ and $h_{125}$ always vary in such a way as to satisfy the equation (\ref{1312341235}), substituting the result of equation (\ref{1312341235}) into $h_{124}$ and $h_{125}$ in $\delta_H(h_{1245})$ gives the following:
\beq
\begin{split}
\delta_H(h_{1234}) &\df \delta_H
\Big(
h_{145}(g_{45}\rhd \big(h_{134}(g_{34}\rhd h_{123})h_{234}^{-1})\big)h_{245}^{-1}h_{235}(g_{35}\rhd h_{123}^{-1})h_{135}^{-1}\Big)
\\
&=   \delta_H
\Big(
\blueunderline{h_{135}^{-1}h_{145}(g_{45}\rhd h_{134})(g_{45}g_{34}\rhd h_{123})}(g_{45}\rhd h_{234}^{-1})h_{245}^{-1}h_{235}(g_{35}\rhd h_{123}^{-1})\Big)\\
&= \delta_H
\Big(
\blueunderline{\redcancel{(g_{35}\rhd h_{123})}\redunderline{h_{135}^{-1}h_{145}(g_{45}\rhd h_{134})}}(g_{45}\rhd h_{234}^{-1})h_{245}^{-1}h_{235}\redcancel{(g_{35}\rhd h_{123}^{-1})}\Big)\\
&=\delta_H
\Big(
\redunderline{h_{345}'}(g_{45}\rhd h_{234}^{-1})h_{245}^{-1}h_{235}\Big) \\
&=\delta_H
\Big(
(h_{345}'(g_{45}\rhd h_{234}^{-1})h_{245}^{-1}h_{235})^{-1}\Big)\\
&=\delta_H
\Big(
h_{235}^{-1}h_{245}(g_{45}\rhd h_{234})h_{345}'^{-1}\Big)\\
&=\delta_H
\Big(
h_{245}(g_{45}\rhd h_{234})h_{345}'^{-1}h_{235}^{-1}\Big)
\end{split}
\label{13r1245}
\eeq
In this transformation, the reason for the red underlined part being able to be simplified is because we previously defined $h_{345}'\df h_{135}^{-1}h_{145}\big(g_{45}\rhd h_{134}\big)$. The reason for being able to simplify as indicated by the blue underlines is that the colors of the 2-simplices $|135|,|145|,|134|$ are correct, i.e. $g_{135}=g_{145}=g_{134}=e$, and this is satisfied through the transformation below.
\beq
\begin{split}
h_{135}^{-1}h_{145}(g_{45}\rhd h_{134})(g_{45}g_{34}\rhd h_{123})
&=h_{135}^{-1}h_{145}(g_{45}\blueunderline{\del(h_{134})g_{34}}\rhd h_{123})(g_{45}\rhd h_{134})\\
&=h_{135}^{-1}h_{145}(g_{45}\blueunderline{g_{14}g_{13}^{-1}}\rhd h_{123})(g_{45}\rhd h_{134})\\
&=h_{135}^{-1}(\redunderline{\del(h_{145})g_{45}g_{14}}g_{13}^{-1}\rhd h_{123})h_{145}(g_{45}\rhd h_{134})\\
&=h_{135}^{-1}(\redunderline{g_{15}}g_{13}^{-1}\rhd h_{123})h_{145}(g_{45}\rhd h_{134})\\
&=(\blueunderline{\del(h_{135}^{-1})g_{15}g_{13}^{-1}}\rhd h_{123})h_{135}^{-1}h_{145}(g_{45}\rhd h_{134})\\
&=(\blueunderline{g_{35}}\rhd h_{123})h_{135}^{-1}h_{145}(g_{45}\rhd h_{134})
\end{split}
\eeq

Using the results from the above equations (\ref{13r145}) and (\ref{13r1245}), equation (\ref{13rhs_2}) can be expressed as follows:

\beq
\begin{split}
Z(r.h.s)
&=|G|^{-5}|H|^{-5}
	  \sum_{g_{13}\in G}
	  \sum_{g_{14}\in G}
	  \sum_{g_{15}\in G}\\
&\quad\times
	 \sum_{h_{123}\in H}
	 \sum_{h_{124}\in H}
	 \sum_{h_{125}\in H}
	 \sum_{h_{134}\in H}
	 \sum_{h_{135}\in H}
	 \sum_{h_{145}\in H}\\
& \quad\quad\times\delta_G(g_{134})\delta_G(g_{135})\delta_G(g_{145})\\
&\quad\quad\quad\times \delta_H(h_{1234})\delta_H(h_{1235})\delta_H(h_{1245})
Z_{remainder}\\
&=|G|^{-5}|H|^{-5}
	  \blueunderline{\sum_{g_{13}\in G}
	 \sum_{h_{123}\in H}
	 \sum_{h_{124}\in H}
	 \sum_{h_{125}\in H}}
	 \redunderline{\sum_{h_{145}\in H}}\\
& \quad\quad\times |G|^2
\delta_G\Big(
\del(h_{345}')g_{45}g_{34}g_{35} ^{-1}
\Big)\\
&\quad\quad\quad\times |H|^2\delta_H\Big(
h_{245}(g_{45}\rhd h_{234})h_{345}'^{-1}h_{235}^{-1}\Big)
Z_{remainder}\\
&=|G|^{-5}|H|^{-5}\blueunderline{ |G||H|^3}
	\redunderline{ \sum_{h_{345}'\in H}}\\
& \quad\quad\times |G|^2
\delta_G\Big(
\del(h_{345}')g_{45}g_{34}g_{35} ^{-1}
\Big)\\
&\quad\quad\quad\times |H|^2\delta_H\Big(
h_{245}(g_{45}\rhd h_{234})h_{345}'^{-1}h_{235}^{-1}\Big)
Z_{remainder}\\
&=|G|^{-2}|H|
\left({1\over |H|}
	{ \sum_{h_{345}'\in H}}
	\right)\\
& \quad\quad\times 
\delta_G\Big(
\del(h_{345}')g_{45}g_{34}g_{35} ^{-1}
\Big)\\
&\quad\quad\quad\times \delta_H\Big(
h_{245}(g_{45}\rhd h_{234})h_{345}'^{-1}h_{235}^{-1}\Big)
Z_{remainder}\\
&=Z(l.h.s)
\end{split}
\label{13rhs_f}
\eeq

As a result, $Z(\text{r.h.s}) = Z(\text{l.h.s})$ has been proved, showing that equation (\ref{ZM}) remains invariant for the 3-dimensional (1-3) move.
\\
\\

\noindent
{\bf (b) 3-dimensional (2-2) move in Figure \ref{2-2mv} }
\\

Finally, for the 3-dimensional (2-2) move in Figure \ref{2-2mv}, let's verify that equation (\ref{ZM}) remains invariant. Similar to before, as the independence on the choice of the total order of 0-simplices has been established, we will perform the calculation based on the total order in Figure \ref{2-2mv}. We designate the left side of Figure \ref{2-2mv} as l.h.s and the right side as r.h.s. In this case, for l.h.s, equation (\ref{ZM}) can be expressed as follows:
\beq
\begin{split}
Z(l.h.s)&=|G|^{-5+9-6}|H|^{5-9+6-2}
\left(
	 {1\over |G|} \sum_{g_{24}\in G}
\right)
\left(
	 {1\over |H|} \sum_{h_{124}\in H}
\right)\\
& \quad\times
\left(
	 {1\over |H|} \sum_{h_{234}\in H}
\right)
\left(
	 {1\over |H|} \sum_{h_{245}\in H}
\right)
\delta_G(g_{124})\delta_G(g_{245})\delta_G(g_{234})\\
& \quad\quad\times \delta_H(h_{1245})\delta_H(h_{1234})
Z_{remainder}
\end{split}
\label{22lhs_1}
\eeq
And on the r.h.s side, it becomes:
\beq
\begin{split}
Z(r.h.s)&=|G|^{-5+9-6}|H|^{5-9+6-2}
\left(
	 {1\over |G|} \sum_{g_{35}\in G}
\right)
\left(
	 {1\over |H|} \sum_{h_{135}\in H}
\right)\\
& \quad\times
\left(
	 {1\over |H|} \sum_{h_{235}\in H}
\right)
\left(
	 {1\over |H|} \sum_{h_{345}\in H}
\right)
\delta_G(g_{135})\delta_G(g_{235})\delta_G(g_{345})\\
& \quad\quad\times \delta_H(h_{1235})\delta_H(h_{1345})
Z_{remainder}
\end{split}
\label{22rhs_1}
\eeq
And here, the term $Z_{remainder}$ is defined as before, representing the common terms between l.h.s and r.h.s.

Let's first simplify the expression for l.h.s. For the 2-simplex $|124|$ in l.h.s, the restriction on $g_{124}$ implies that $g_{124} = e$ when the color is correct. Therefore, let $g_{24}$ vary in such a way that it satisfies the following:
\beq
\begin{split}
g_{124} &\df \del(h_{124})g_{24}g_{12}g_{14}^{-1} = e\\
&\iff g_{24}=\del(h_{124}^{-1})g_{14}g_{12}^{-1}
\end{split}
\label{22g24_124}
\eeq
Using the expression (\ref{22g24_124}), by substituting it into the restrictions $g_{234}, g_{245}$ for the 2-simplices $|234|$ and $|245|$, we obtain the following:
\beq
\begin{split}
g_{234} &\df \del(h_{234})g_{34}g_{23}g_{24}^{-1}\\
&=\del(h_{234})g_{34}g_{23}g_{12}g_{14}^{-1}\del(h_{124})
\end{split}
\label{22g234_1}
\eeq
\beq
\begin{split}
g_{245} &\df \del(h_{245})g_{45}g_{24}g_{25}^{-1}\\
&= \del(h_{245})g_{45}\del(h_{124}^{-1})g_{14}g_{12}^{-1}g_{25}^{-1}
\end{split}
\label{22g245_1}
\eeq
Here, for the restrictions $h_{1245}, h_{1234}$ of the 3-simplex $|1245|$ and $|1234|$ in l.h.s, when the colors are correct, $h_{1245} = h_{1234} = e$ holds, and therefore, the following equation is satisfied:
\beq
\begin{split}
h_{1234} &\df h_{134}(g_{34}\rhd h_{123})h_{234}^{-1}h_{124}^{-1} = e \\
&\iff h_{234} =h_{124}^{-1}h_{134}(g_{34}\rhd h_{123})
\label{22h1234_234}
\end{split}
\eeq
\beq
\begin{split}
h_{1245} &\df h_{145}(g_{45}\rhd h_{124})h_{245}^{-1}h_{125}^{-1} = e \\
&\iff h_{245} =h_{124}^{-1}h_{145}(g_{45}\rhd h_{125})
\end{split}
\label{22h1245_245}
\eeq
Here, substituting the expression (\ref{22h1234_234}) into the expression (\ref{22g234_1}) for $h_{234}$, we get the following:
\beq
\begin{split}
g_{234} &=\del(h_{234})g_{34}g_{23}g_{12}g_{14}^{-1}\del(h_{124})\\
&=\del\Big(h_{124}^{-1}h_{134}(g_{34}\rhd h_{123})\Big)g_{34}g_{23}g_{12}g_{14}^{-1}\del(h_{124})\\
&=\del(h_{124}^{-1})\del(h_{134})g_{34}\blueunderline{\del(h_{123})g_{23}g_{12}}g_{14}^{-1}\del(h_{124})\\
&=\del(h_{124}^{-1})\redunderline{\del(h_{134})g_{34}\blueunderline{g_{13}}}g_{14}^{-1}\del(h_{124})\\
&=\del(h_{124}^{-1})\redunderline{g_{14}}g_{14}^{-1}\del(h_{124})\\
&=e
\end{split}
\label{22g234_2}
\eeq
When transforming the expression (\ref{22g234_2}), the correctness of the color factors for the 2-quartets $|123|$ and $|134|$ was ensured, namely, by using $g_{123}=g_{134}=e$. Substituting (\ref{22h1245_245}) into (\ref{22g245_1}) for $h_{245}$, we get the following:
\beq
\begin{split}
g_{245} &= \del(h_{245})g_{45}\del(h_{124}^{-1})g_{14}g_{12}^{-1}g_{25}^{-1}\\
&=\del\Big(h_{125}^{-1}h_{145}(g_{45}\rhd h_{124})\Big)g_{45}\del(h_{124}^{-1})g_{14}g_{12}^{-1}g_{25}^{-1}\\
&=\del(h_{125}^{-1})\blueunderline{\del(h_{145})g_{45}}\redcancel{\del(h_{124})\del(h_{124}^{-1})}\blueunderline{g_{14}}g_{12}^{-1}g_{25}^{-1}\\
&=\redunderline{\del(h_{125}^{-1})\blueunderline{g_{15}}g_{12}^{-1}}g_{25}^{-1}\\
&=\redunderline{g_{25}}g_{25}^{-1}\\
&=e
\end{split}
\label{22g245_2}
\eeq
Using $g_{145} = g_{125} = e$ in the transformation of (\ref{22g245_2}), we get the following:

Now, let's consider the transformation of $h_{1245}$ and $h_{1234}$ in l.h.s. Since $h_{1234}=e$ due to the correct coloring of the 3-simplex $|1234|$, let's assume that $h_{124}$ moves in such a way that it always satisfies the following:
\beq
\begin{split}
h_{1234} &\df h_{134}(g_{34}\rhd h_{123})h_{234}^{-1}h_{124}^{-1} = e \\
&\iff  h_{124}=h_{134}(g_{34}\rhd h_{123})h_{234}^{-1}
\end{split}
\label{22h124_1234}
\eeq

Using the expression (\ref{22h124_1234}), substituting into $h_{1245}$ for $h_{124}$, we get the following:
\beq
\begin{split}
h_{1245} &\df h_{145}(g_{45}\rhd h_{124})h_{245}^{-1}h_{125}^{-1}  \\
&=h_{145}
\left(
g_{45}\rhd
	\left(
	h_{134}(g_{34}\rhd h_{123})h_{234}^{-1}
	\right)
\right)h_{245}^{-1}h_{125}^{-1}\\
&=h_{145}
\left(
g_{45}\rhd h_{134}
\right)
\left(
g_{45}g_{34}\rhd h_{123}
\right)
\redunderline{\left(
g_{45}\rhd h_{234}^{-1}
\right)
h_{245}^{-1}}h_{125}^{-1}\\
&=h_{145}
\left(
g_{45}\rhd h_{134}
\right)
\left(
g_{45}g_{34}\rhd h_{123}
\right)
\redunderline{h^{-1}}h_{125}^{-1}
\end{split}
\label{22h1245}
\eeq

However, in the last equation, we express it as $h\df h_{245}\left(g_{45}\rhd h_{234}\right)$. Using the results from equations (\ref{22g24_124}) to (\ref{22h1245}), the expression (\ref{22lhs_1}) can be written as follows:
\beq
\begin{split}
Z(l.h.s)&=|G|^{-2}
\blueunderline{\left(
	 {1\over |G|} \sum_{g_{24}\in G}
\right)
\left(
	 {1\over |H|} \sum_{h_{124}\in H}
\right)}
\redunderline{\left(
	 {1\over |H|} \sum_{h_{234}\in H}
\right)
\left(
	 {1\over |H|} \sum_{h_{245}\in H}
\right)}\\
& \quad\times
\delta_G(g_{124})\delta_G(g_{245})\delta_G(g_{234}) \delta_H(h_{1245})\delta_H(h_{1234})
Z_{remainder}\\
&=|G|^{-2}
\blueunderline{	 {1\over |G|} 
	 {1\over |H|}}
\redunderline{\left(
	 {1\over |H|} \sum_{h\in H}
\right)
}\\
& \quad\times
|G|^3|H|\delta_H(h_{145}
\left(
g_{45}\rhd h_{134}
\right)
\left(
g_{45}g_{34}\rhd h_{123}
\right)
h^{-1}h_{125}^{-1})
Z_{remainder}\\
\end{split}
\label{22lhs_2}
\eeq
\\

Next, similarly, let's transform the r.h.s. When it comes to the color factor $g_{135}$ for the 2-simplex $|135|$ on the r.h.s., if its color is correct, i.e., $g_{135}=e$, we assume that $g_{35}$ adjusts itself to satisfy the following:
\beq
\begin{split}
g_{135} &\df \del(h_{135})g_{35}g_{13}g_{15}^{-1} = e\\
&\iff g_{35}=\del(h_{135}^{-1})g_{15}g_{13}^{-1}
\end{split}
\label{22g35_135}
\eeq
Using the expression (\ref{22g35_135}), substituting into $g_{35}$ in $g_{235}$ and $g_{345}$ we get next formula:
\beq
\begin{split}
g_{235} &\df \del(h_{235})g_{35}g_{23}g_{25}^{-1}\\
& \iff\del(h_{235})\del(h_{135}^{-1})g_{15}g_{13}^{-1}g_{23}g_{25}^{-1}
\end{split}
\label{22g235_1}
\eeq
\beq
\begin{split}
g_{345} &\df \del(h_{345})g_{45}g_{34}g_{35}^{-1}\\
&= \del(h_{345})g_{45}g_{34}g_{13}g_{15}^{-1}\del(h_{135})
\end{split}
\label{22g345_1}
\eeq

Here, for the l.h.s. expressions $h_{1235}$ and $h_{1345}$, when the color factors are correct, i.e., $h_{1235}=h_{1345}=e$, the following equation holds.
\beq
\begin{split}
h_{1235} &\df h_{135}(g_{35}\rhd h_{123})h_{235}^{-1}h_{125}^{-1} = e \\
&\iff h_{235}= h_{125}^{-1}h_{135}(g_{35}\rhd h_{123})
\label{22h1235_235}
\end{split}
\eeq
\beq
\begin{split}
h_{1345} &\df h_{145}(g_{45}\rhd h_{134})h_{345}^{-1}h_{135}^{-1} = e \\
&\iff h_{345}=h_{135}^{-1}h_{145}(g_{45}\rhd h_{134})
\end{split}
\label{22h1345_345}
\eeq
Here, substituting the expression (\ref{22h1235_235}) into the equation for $h_{235}$ in (\ref{22g235_1}), we obtain the following:
\beq
\begin{split}
g_{235} &= \del\left(h_{235}\right)\del(h_{135}^{-1})g_{15}g_{13}^{-1}g_{23}g_{25}^{-1}\\
&=\del\left(h_{125}^{-1}h_{135}(g_{35}\rhd h_{123})\right)\del(h_{135}^{-1})g_{15}g_{13}^{-1}g_{23}g_{25}^{-1}\\
&=\del(h_{125}^{-1})\redcancel{\del(h_{135})\del(h_{135}^{-1})}g_{15}\blueunderline{g_{13}^{-1}\del(h_{123})g_{23}}g_{25}^{-1}\\
&=\redunderline{\del(h_{125}^{-1})g_{15}\blueunderline{g_{12}^{-1}}}g_{25}^{-1}\\
&=\redunderline{g_{25}}g_{25}^{-1}\\
&=e
\end{split}
\label{22g235_2}
\eeq
When transforming the expression (\ref{22g235_2}), the color factors $g_{123}=g_{125}=e$ were used. Similarly, substituting the expression (\ref{22h1345_345}) into $h_{345}$ in equation (\ref{22g345_1}) results in the following:
\beq
\begin{split}
g_{345} &= \del\left(h_{345}\right)g_{45}g_{34}g_{13}g_{15}^{-1}\del(h_{135})\\
&=\del\left(h_{135}^{-1}h_{145}(g_{45}\rhd h_{134})\right)g_{45}g_{34}g_{13}g_{15}^{-1}\del(h_{135})\\
&=\del(h_{135}^{-1})\del(h_{145})g_{45}\redunderline{\del(h_{134})g_{34}g_{13}}g_{15}^{-1}\del(h_{135})\\
&=\del(h_{135}^{-1})\blueunderline{\del(h_{145})g_{45}\redunderline{g_{14}}}g_{15}^{-1}\del(h_{135})\\
&=\del(h_{135}^{-1})\blueunderline{g_{15}}g_{15}^{-1}\del(h_{135})\\
&=e
\end{split}
\label{22g345_2}
\eeq

When transforming the expression (\ref{22g245_2}), the color factors $g_{134}=g_{145}=e$ were used.

Next, let's consider the transformation of the r.h.s. expressions $h_{1235}$ and $h_{1345}$. Due to the correct color factor for the 3-simplex $|1345|$, i.e., $h_{1345}=e$, we use this to ensure that $h_{135}$ always satisfies the following:
\beq
\begin{split}
h_{1345} &\df h_{145}(g_{45}\rhd h_{134})h_{345}^{-1}h_{135}^{-1} = e \\
&\iff  h_{135}=h_{145}(g_{45}\rhd h_{134})h_{345}^{-1}
\end{split}
\label{22h135_1345}
\eeq

Using the expression (\ref{22h135_1345}), substituting into $h_{1235}$ for $h_{135}$ yields the following:
\beq
\begin{split}
h_{1235} &\df h_{135}(g_{35}\rhd h_{123})h_{235}^{-1}h_{125}^{-1}  \\
&=h_{145}(g_{45}\rhd h_{134})\blueunderline{h_{345}^{-1}(g_{35}\rhd h_{123})}h_{235}^{-1}h_{125}^{-1} \\
&=h_{145}(g_{45}\rhd h_{134})\blueunderline{\left(\redunderline{\del(h_{345}^{-1})g_{35}}\rhd h_{123}\right)h_{345}^{-1}}h_{235}^{-1}h_{125}^{-1} \\
&=h_{145}(g_{45}\rhd h_{134})\left(\redunderline{g_{45}g_{34}}\rhd h_{123}\right)h_{345}^{-1}h_{235}^{-1}h_{125}^{-1}\\
&= h_{145}(g_{45}\rhd h_{134})\left(g_{45}g_{34}\rhd h_{123}\right)h'^{-1}h_{125}^{-1}
\end{split}
\label{22h1235}
\eeq
However, in the last equation, we express it as $h'\df h_{235}h_{345}$. In the underlined portion in red, due to the correct color factor for the 2-simplex $|345|$, i.e., $g_{345}=e$, this holds.

Using the results from equations (\ref{22g35_135}) to (\ref{22h1235}), the expression (\ref{22rhs_1}) can be written as follows:
\beq
\begin{split}
Z(r.h.s)&=|G|^{-2}
\blueunderline{\left(
	 {1\over |G|} \sum_{g_{35}\in G}
\right)
\left(
	 {1\over |H|} \sum_{h_{135}\in H}
\right)}
\redunderline{\left(
	 {1\over |H|} \sum_{h_{235}\in H}
\right)
\left(
	 {1\over |H|} \sum_{h_{345}\in H}
\right)}\\
& \quad\times
\delta_G(g_{135})\delta_G(g_{235})\delta_G(g_{345}) \delta_H(h_{1235})\delta_H(h_{1345})
Z_{remainder}\\
&=|G|^{-2}
\blueunderline{	 {1\over |G|} 
	 {1\over |H|}}
\redunderline{\left(
	 {1\over |H|} \sum_{h'\in H}
\right)
}\\
& \quad\times
|G|^3|H|\delta_H\left(h_{145}(g_{45}\rhd h_{134})\left(g_{45}g_{34}\rhd h_{123}\right)h'^{-1}h_{125}^{-1}\right)
Z_{remainder}
\end{split}
\label{22rhs_2}
\eeq
\\

Comparing equation (\ref{22lhs_2}) with equation (\ref{22rhs_2}), 
it is evident that the forms of the equations are identical. Therefore, 
we can conclude that $Z(\text{l.h.s.})=Z(\text{r.h.s.})$. 
Thus, we have established the invariance of equation (\ref{ZM}) 
with respect to the moves depicted in Figure \ref{2-2mv}.

\section*{Acknowledgement}

I would like to express my deep gratitude for my supervisor Professor Yuji Terashima. He gave me a really large amount of worthful advices and suggestions,
and offered a very comfortable place that made my life of research so pleasant.
Without his continual help and encouragement, this paper would not be completed.

\end{document}